\documentclass[12pt]{article}
\usepackage{array,CJK}
\usepackage{latexsym,epsfig}
\usepackage{amssymb,amsmath,amsthm}
\usepackage[utf8]{inputenc}
\usepackage{graphicx} 
\usepackage{mwe}
\usepackage{subfig}
\usepackage{authblk}
\newcommand{\al}{\alpha}

\newcommand{\be}{\beta}
\newcommand{\la}{\lambda}
\newcommand{\ep}{\epsilon}

\newcommand{\vp}{\varphi}
\newcommand{\sig}{\sigma}
\renewcommand{\th}{\theta}
\newcommand{\Th}{\Theta}
\newcommand{\rme}{{\rm e}}
\newcommand{\bfN}{\mathbb{N}}
\newcommand{\bfR}{\mathbb{R}}

\newcommand{\bfZ}{\mathbb{Z}}
\newcommand{\bfC}{\mathbb{C}}
\newcommand{\bk}{{\bf k}}
\newcommand{\bp}{{\bf p}}
\newcommand{\bq}{{\bf q}}

\newcommand{\bv}{{\bf v}}
\newcommand{\by}{{\bf y}}

\usepackage{color}

\title{Dispersion relations of periodic quantum graphs associated with Archimedean tilings (I)}
\author[1]{Yu-Chen Luo}
\author[1,2] {Eduardo O.\ Jatulan}
\author[1] {Chun-Kong Law}
 \affil[1]{ \footnotesize{Department of Applied Mathematics, National Sun Yat-sen University,
Kaohsiung, Taiwan 80424. Email: law@math.nsysu.edu.tw}}
\affil[2]{ \footnotesize{Institute of Mathematical Sciences and Physics, University of the Philippines Los Banos, Philippines  4031.
Email: eojatulan@up.edu.ph}}

\begin{document}
\maketitle
 \begin{abstract}
  There are totally 11 kinds of Archimedean tiling for the plane. Applying the Floquet-Bloch theory, we derive the dispersion relations of the  periodic quantum graphs associated with a number of Archimedean tiling, namely the triangular tiling  {$(3^6)$}, the elongated triangular tiling  {$(3^3,4^2)$}, the trihexagonal tiling  {$(3,6,3,6)$} and the truncated square tiling  {$(4,8^2)$}. The derivation makes use of characteristic functions, with the help of the symbolic software Mathematica.
 The resulting dispersion relations are surprisingly simple and symmetric.  They show that in each case the spectrum is composed of point spectrum and an absolutely continuous spectrum.  We further analyzed on the structure of the absolutely continuous spectra.  Our work is motivated by the studies on the periodic quantum graphs associated with hexagonal tiling  in \cite{KP} and \cite{KL}.
\\[0.1in]
 \end{abstract}
 \vskip1in
 \newpage
 \section{Introduction}
 \setcounter{equation}{0}\hskip0.25in
 Recently there have been a lot of studies on quantum graphs, which is essentially the spectral problem of a one-dimensional
 Schr\"{o}dinger operator acting on the edge of a graph, while the functions have to satisfy some boundary conditions as well
 as vertex conditions which are usually the continuity and Kirchhoff conditions. Quantum graphs finds its applications in nanomaterials,
 network theory and chemistry.  {Interested readers may consult  \cite{EKKST,BK} for a broad introduction to
 quantum graphs.} In particular, we developed a reduction formula for the characteristic functions \cite{LP, LY} whose
 zeros are exactly the square roots of eigenvalues for compact quantum trees. This characteristic function approach seems to be useful in the study of spectra for quantum graphs.
 In this study, we shall employ it to derive the dispersion relations of periodic quantum graphs associated with Archimedean tiling \cite{GS}.

 This paper is motivated by the work of Kuchment-Post \cite{KP}.
 On the plane, there are plenty of ways to form a tessellation with different patterns \cite{GS}. A tessellation is formed
 by regular polygons having the same edge length such that these regular polygons will fit around at
 each vertex and the pattern at every vertex is isomorphic, then it is  called a \textit{uniform tiling} (also called  Archimedean tilings in literature).
 There are 11 types of uniform tilings,  as shown in Table \ref{Tab5.1}.
  We try to investigate the periodic spectrum of some differential operators, in particular the Schr\"{o}dinger operator acting on these infinite
  graphs. Among these 11 types, the square tiling
 (with vertex configuration ($4^4$)) was studied in \cite{LLW}.  The hexagonal tiling was solved by Kuchment-Post \cite{KP} and Korotyaev-Lobanov \cite{KL}
 independently.  Using Floquet-Bloch theory, we  {have} derived and analyzed on the dispersion relations of the remaining periodic quantum graphs in a systematic way.
 The analytic variety of these relations, also called Bloch variety, gives the spectrum of the differential
 operators.

 Given an infinite graph $G=E(G)\cup V(G)$ generated by an Archimedean tiling, we let $H$ denote a Schr\"odinger operator on , i.e.,
 $$
 H \, \by(x)=-\frac{d^2}{d x^2}\by(x)+\bq(x)\, \by(x),
 $$
  {where $ \bq\in L^2_{loc}(G)$ is periodic on the tiling (explained below), and} the domain $D(H)$ consists of all  {admissible} functions $\by(x)$ (union of functions $y_e$ for each edge $e\in E(G)$) on $G$ in the sense that
 \begin{enumerate}
 \item[(i)] $\displaystyle   y_e\in {\cal H}^2(e)$ for all $e\in E(G)$;
 \item[(ii)] $\displaystyle   \sum_{ {e\in E(G)}} \| y_e\|^2_{{\cal H}^2(e)} <\infty$;
 \item[(iii)] Neumann vertex conditions (or continuity-Kirchhoff conditions at vertices), i.e., for any vertex $\bv\in V(G)$,
 $$
 y_{e_1}(\bv)=y_{e_2}(\bv)\qquad \mbox{ and }\qquad \sum_{ {e\in E(G)}} y_e'(\bv)=0.
 $$
 Here $y_e'$ denotes the directional derivative along the edge $e$ from $\bv$,  {and $\| \cdot\|^2_{{\cal H}^2(e)}$ denotes the Sobolev norm of 2 distribution derivatives.}
 \end{enumerate}
 The defined operator $H$ is well known to be unbounded and self-adjoint. As $H$ is also periodic, we may apply
 the Floquet-Bloch theory \cite{BK,E73,RS} in the study of its spectrum.
  {Let $\vec{k}_1,\vec{k}_2$ be two linear independent vectors in $\bfR^2$. Define $\vec{\bk}=(\vec{k}_1,\vec{k}_2)$, and $\bp=(p_1,p_2)\in \bfZ^2$.
 For any set $S\subset G$, we define $\bp \circ S$ to be an action on $S$ by a shift of $\bp\cdot\bk=p_1\bk_1+p_2\bk_2$.
 A compact set $W\subset G$ is said to be a fundamental domain if \
 $$
 G=\bigcup\{ \bp\circ W:\ \bp\in \bfZ^2\},
 $$
 and for any different $\bp,\bp'\in \bfZ^2$, $(\bp\circ W)\cap (\bp'\circ W)$ is a finite set in $G$. The potential function $\bq$ is said to be periodic if
 $\bq(x+\bp\cdot \vec{\bk})=\bq(x)$ for all $\bp\in \bfZ^2$ and all $x\in G$.}
  Take the quasi-momentum $\Theta=(\theta_1,\theta_2)$ in the Brillouin zone
 $B=[-\pi,\pi]^2$. Let $H^{\Th}$ be the Bloch  {Hamiltonian} that acts on  {$L^2(G)$, and the dense domain $D(H^{\Th})$ consists of admissible functions $\by$ which satisfy the Floquet-Bloch condition}
 \begin{equation}
 \by(x+\bp\cdot \vec{\bk})= \rme^{i(\bp\cdot\Theta)}\by(x),
 \label{eq1.02}
 \end{equation}
 for all $\bp\in \bfZ^2$ and all $x\in G$.   {Such functions are uniquely determined by their restrictions on the fundamental domain $W$. Hence for fixed $\Th$,}
 the operator $H^\Th$ has purely discrete spectrum $\sig(H^\Th)=\{\la_j(\Th):\ j\in \bfN\}$, where
 $$
  {\la_1(\Th)\leq \la_2(\Th)\leq \cdots \leq \la_j(\Th)\leq\cdots},\quad \mbox{ and } \la_j(\Th)\rightarrow\infty\mbox{ as } j\to\infty.
 $$

 We say $H$ is a direct integral of $\displaystyle   H^\Th$, denoted as
 $$
 H=\int_B^\oplus H^{\Th}\, d\Th.
 $$
 Moreover, by \cite{BK,RS},
 $$
 \sig(H)= \bigcup \{\sig(H^{\Th}):\ \Th\in [-\pi,\pi]^2\}.
 $$
 Furthermore,  {it is known that  singular continuous spectrum is absent in $\sig(H)$. So  $\sig(H)$ consists of only point spectrum and absolutely continuous spectrum \cite[Theorem 4.5.9]{K93}}.

 On the interval $[0,a]$, we let $C(x,\rho)$ and $S(x,\rho)$  {($\la=\rho^2$)} be the solutions of
 $$
 -y''+q y=\lambda y
 $$
 such that $C(0,\rho)=S'(0,\rho)=1$, $C'(0,\rho)=S(0,\rho)=0$. We call them cosine-like function and sine-like function
 respectively. In particular,
 \begin{eqnarray*}
 C(x,\rho)&=&\displaystyle\cos(\rho x)+\frac1{\rho}\int^x_0\sin(\rho (x-t))q(t)C(t, \rho)dt,\\
 S(x, \rho)&=&\displaystyle \frac{\sin(\rho x)}{\rho}+\frac1{\rho}\int^x_0\sin(\rho (x-t))q(t)S(t, \rho)dt.
 \end{eqnarray*}
  So for the periodic quantum graphs above, suppose the edges $\{ e_1,\ldots,e_I\}$ lie in a typical fundamental domain $W$, while
  $(q_1,\ldots,q_I)$ are the potential functions acting on these edges.
 Assume also that the potential functions $q_i$'s are identical and even.
 Then the dispersion relation, which defines the spectral value $\la$ as a function of the quasimomentum $\Th=(\theta_1,\theta_2)$, of the periodic quantum graph associated with square tiling was found \cite{LLW} to be
 $$
 S(a,\rho)^2\, (S'(a,\rho)^2-\cos^2(\frac{\theta_1}{2})\cos^2(\frac{\theta_2}{2}))=0.
 $$
 The dispersion relation for hexagonal tiling \cite{KP,KL} is
 $$
 S(a,\rho)^2\,\left( 9 S'(a,\rho)^2-1-8\cos(\frac{\theta_1}{2})\cos(\frac{\theta_2}{2})\cos(\frac{\theta_1-\theta_2}{2})\right)=0.
 $$
 Thus the spectra of the two periodic quantum graphs are expressed in terms of the sine-like functions defined on the line.
 The function $S'(a,\rho)$ in the above dispersion relation is exactly the function $\displaystyle \eta(\la)=\frac{1}{2}D(H)$ (the Hill's discriminant)
 defined in \cite{KP,KL}.

 In this study, we shall use a characteristic function approach to derive the dispersion relations of some of the other periodic quantum graphs associated with Archimedean tilings, namely the triangular tiling ($3^6$),
 the elongated triangular tiling ($3^3,\, 4^2$), the truncated square tiling $(4,8^2)$ and the trihexagonal tiling $(3,6,3,6)$. These tilings are denoted in short forms  as
 $T$, $eT$, $trS$, and $TH$ respectively.  {We shall also assume that the periodic function $\bq$ may have components $q_i$ defined on each edge $e_i$ which might neither be identical, nor even.} In this way,
 we obtain more general dispersion relations.  As
 the computation is more complicated, we need the help of the symbolic software Mathematica. For example, the dispersion relation for the triangular tiling
 is
 \begin{eqnarray*}
 \lefteqn{ S_1'S_2S_3+S_1 S_2'S_3+S_1 S_2 S_3'+C_1 S_2S_3+C_2S_1S_3+C_3S_1S_2}\\
  &&-2S_1 S_2\, \cos(\theta_1-\theta_2)-2 S_2 S_3\, \cos\theta_1-2 S_1 S_3\cos\theta_2\  = \ 0.
  \end{eqnarray*}
  \noindent
  (cf.\ Theorem \ref{th2.1}). From now on, we shall use the following abbreviated symbols,
  $$
 S_j=S_j(a,\rho),\ S_j'=S_j'(a,\rho),\ C_j=C_j(a,\rho),\ C_j'=C_j'(a,\rho);
 $$
  and similar symbols for $S,\ S',\ C,\ C'$. If we assume that potential functions are identical and even, then we arrive at
  a simple and symmetric dispersion relation as below. The following theorem is our main result.
 \newtheorem{th1.0}{Theorem}[section]
 \begin{th1.0}
 \label{th1.0}
  Assume that all the $q_j$'s are identical (denoted as $q$), and even. We also let $\theta_1,\ \theta_2\in[-\pi,\pi]$.
 \begin{enumerate}
 \item[(a)] For the triangular tiling  {$(3^6)$}, the dispersion relation of the associated periodic quantum graph is
 \begin{equation}
 S^2 \left(3S'+1-4 \cos(\frac{\theta_1}{2})\cos(\frac{\theta_2}{2})\cos(\frac{\theta_2-\theta_1}{2})\right)=0.
 \label{eq1.03}
 \end{equation}
 \item[(b)] For the elongated triangular tiling  {$(3^3,4^2)$}, the dispersion relation of the associated periodic quantum graph is given by
  \begin{equation}
  S^3\, \{25(S')^2-20\cos\theta_1S'-8\cos(\frac{\theta_1}{2})\cos(\frac{\theta_2}{2})\cos(\frac{\theta_1-\theta_2}{2})+4\cos^2\theta_1-1\}=0.
  \label{eq1.04}
  \end{equation}
 \item[(c)] For the truncated square tiling  {$(4,8^2)$}, the dispersion relation of the associated periodic quantum graph is given by
  \begin{equation}
 S^2\left\{ 81 S'^4-54 S'^2-12S'(\cos\theta_1+\cos\theta_2)+1-4\cos\theta_1\, \cos\theta_2\right\} = 0.
  \label{eq1.05}
  \end{equation}
 \item[(d)] For trihexagonal tiling  {$(3,6,3,6)$}, the dispersion relation of the associated periodic quantum graph is given by
  \begin{equation}
 S^3(2S'+1)\, \left(2 S'^2-S'-\cos\frac{\theta_2}{2}\cos\frac{\theta_1-\theta_2}{2})\right)=0
  \label{eq1.06}
  \end{equation}
 \end{enumerate}
 \end{th1.0}
 Since these dispersion relations are explicit and simple, it is possible to do  further analysis to understand the spectra.
 They show that in each case the spectrum is composed of point spectrum ($\sig_p$) which have an infinite number
 of eigenfunctions, and an absolutely continuous spectrum ($\sig_{ac}$).  In fact, eigenfunctions associated with point spectrum can be easily constructed through the functions
 $S=S(a,\rho)=0$, or $S'=S'(a,\rho)=0$. Also all the absolutely continuous spectra obviously have a band and gap structure.

 The Archimedean tilings provide good models for crystal lattices.
 It has been known that molecules like graphene
 or boron nitride (BN) have crystal lattices \cite{ALM, LAM}. In fact in graphene the carbon atoms are located at vertices of regular hexagons. By a
 quantum network model (QNM), the wave functions at the bonds (Bloch waves) satisfy a Schr\"{o}dinger equation along the lines, and
 continuity and Kirchhoff conditions at the vertices. So under this QNM, the spectrum represent the energy of the wave functions in
 graphene. Thus the associated spectral analysis has physical significance in quantum mechanics. In fact, in \cite{ALM}, the potential function for graphene is given as
 $$
 q(x)=-0.85+\displaystyle\frac{d}{1.34}\sin^2\left(\frac{\pi x}{d}\right)
 $$
 where $d$ is the distance between neighboring carbon atoms, and $d=1.43$ {\AA}.  Hence the spectra of their periodic graphs, as in the case of graphene, have physical meaning and their analyses worthwhile. \par

 Graphene is famous for its unusual electric and mechanical properties \cite{G2011,N2011}.  In \cite{EI}, 14 different 2D carbon allotropes were proposed to be likely to
 form novel graphene-like materials. These allotropes include two lattices from the Archimedean tiling - truncated square tiling $(4,8,8)$ and
 truncated trihexagonal tiling $(4,6,12)$. In  \cite{K2013}, Do and Kuchment studied another allotrope called graphyne (similar to $(4,6,4,6)$ but
 with  rhombi instead of squares in between rows of hexagons). They derived the dispersion relation of this graphyne and analyzed on the spectrum.  Their method also works for
 our models.

 \newtheorem{th1.1}[th1.0]{Theorem}
 \begin{th1.1}
 \label{th1.1}
 Assuming all $q_i's$ are identical and even,
 \begin{enumerate}
 \item[(a)] $\displaystyle  \sigma_{ac}(H_{T})=\left\{\rho^2\in \bfR:\ S'(a,\rho)\in\left[-\frac{1}{2},1\right]\right\}=\ \sigma_{ac}(H_{TH})$
 \item[(b)] $\displaystyle  \sigma_{ac}(H_{eT})=\left\{\rho^2\in \bfR:\ S'(a,\rho)\in\left[-\frac{3}{5},1\right]\right\}$
 \item[(c)]   $\displaystyle   \sigma_{ac}(H_{trS})=\left\{\rho^2\in \bfR:\ S'(a,\rho)\in\left[-1,1\right]\right\}$
 \end{enumerate}
\end{th1.1}

  It is easy to see that
 $$
 \sig_{ac}(H_S)=\{\rho^2\in \bfR:\ S'(a,\rho)\in [-1,1]\} =\sig_{ac}(H_{hex}).
 $$
 It means that the  {absolutely continuous spectra} for these two Archimedean tilings are exactly the same as those for the Hill operator.
 It was shown in \cite{S76} for almost all $C^\infty$ periodic potential $q_0$ on $\bfR$, the spectrum for
 the Hill operator opens up at each spectral value at the edges $S'(a,\rho)=\pm 1$, upon a perturbation of the potential.
 Furthermore, Do-Kuchment \cite{K2013} proved that for graphyne, in the free case when $q=0$, there is a conical singularity at these special points, which are called
 Dirac points. It was these Dirac points that account for the special electronic properties of graphene.   {Fefferman and Weinstein \cite{FW12} studied the existence of
 infinitely many Dirac points for the hexagonal tiling. In \cite{BC}, Berkolaiko and Comech  associated the existence of Dirac points with Berry phases.}
 We found that
 for all the above quantum graphs  there are points such that $S'(a,\rho)=1$, and for truncated square tiling, it is possible to have $S'(a,\rho)=-1$.  We shall study this
 Dirac point issue in more detail in another paper.

 In the following sections, we shall study the characteristic functions and dispersion relations for the periodic quantum graphs associated with triangular tiling,
 elongated triangular tiling, truncated square tiling and trihexagonal tiling. These derivations, which constitutes the proof of Theorem~\ref{th1.1}, will
 be given for each tiling from sections 2 through section 5.
 Then in section 6, we shall study the spectra thus obtained in more detail. There are five Archimedean tilings left.  We shall deal with them in another paper.

  It is interesting to note that many of these dispersion relations involves the expression $\cos a+\cos b+\cos (a-b)$.
 It is easy to see that
  \begin{eqnarray}
  \left|1+\rme^{i\theta_1}+\rme^{i\theta_2}\right|^2&=& 3+2(\cos \theta_1+\cos \theta_2+\cos (\theta_1-\theta_2))\nonumber\\
  &=&1+ 8 \cos (\frac{\theta_1}{2}) \cos (\frac{\theta_2}{2})\cos(\frac{\theta_1-\theta_2}{2}). \label{eq1.2}
  \end{eqnarray}

  Assume that $C$ and $S$ are solutions of the Sturm-Liouville equation on the interval $(0,1)$ with initial conditions
  $$
  C(0,\rho)=1=S'(0,\rho),\qquad S(0,\rho)=C'(0,\rho)=0,
  $$
  where $\la=\rho^2\in \bfC$.  Part (d) below is extremely useful in this paper.
  \newtheorem{th1.2}[th1.0]{Theorem}
  \begin{th1.2}
  \label{th1.2}
  Suppose that the potential function $q$ is even. Then for any $\la=\rho^2$ in $\bfC$,
  \begin{enumerate}
  \item[(a)] $\displaystyle  C(a,\rho)=2C(\frac{a}{2},\rho)S'(\frac{a}{2},\rho)-1=1+2\, S(\frac{a}{2},\rho)C'(\frac{a}{2},\rho)$;
  \item[(b)] $\displaystyle  S(a,\rho)= 2\, S(\frac{a}{2},\rho)S'(\frac{a}{2},\rho)$;
  \item[(c)] $\displaystyle  C'(a,\rho) = 2\, C(\frac{a}{2},\rho)S'(\frac{a}{2},\rho)$;
  \item[(d)] $ \displaystyle  S'(a,\rho)=C(a,\rho)$.
  \end{enumerate}
  \end{th1.2}
  The above theorem was given in Chapter 1 of classical monograph of Magnus and Winkler \cite[p.8]{MW}.
  Recently it was reproved by Pivovarchik-Rozhenko
   \cite{PR13} using theory of entire functions.

   Besides the abbreviation of $S,S',C,C'$, we introduce some more notations to simplify the long expressions. For example
   \begin{equation}
   \left\{ \begin{array}{rcl}
   (CSS)_{ijk} &:=& C_iS_jS_k+C_jS_kS_i+C_kS_iS_j \\
   (CS)_{ij'} &:=& C_i S_j'+S_iC_j'
   \end{array} \right.\ .
   \label{eq1.3}
   \end{equation}
   These notations will show up later.
 \section{Triangular tiling}
\begin{figure}[h!]
 \centering\includegraphics[width=1\textwidth]{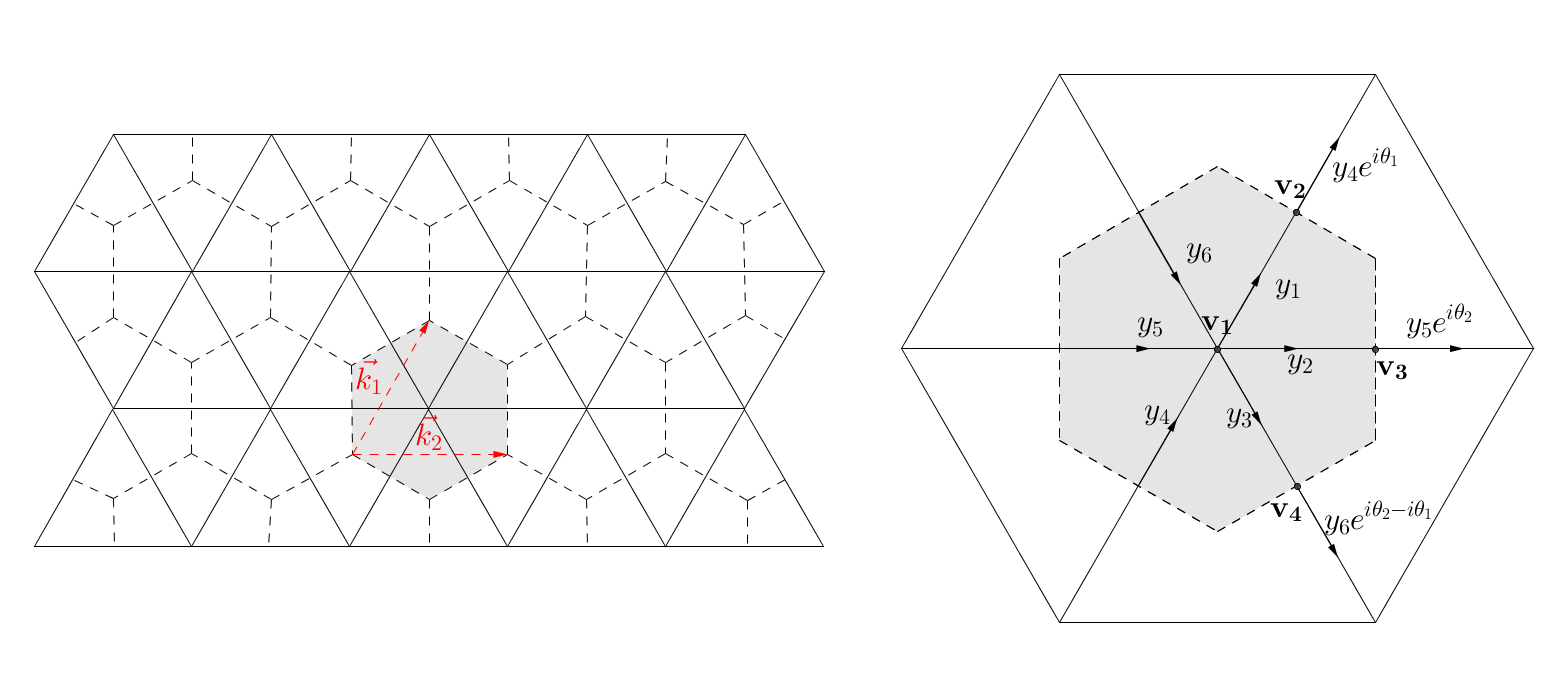}
 \caption{Fundamental domain for triangular tiling}
 \label{fig4.3}
 \end{figure}

 \setcounter{equation}{0}\hskip0.25in
  The fundamental domain for triangular tiling is quite different from that of hexagonal tiling.
  Its edges cut through the triangular edges.  We let it be a regular hexagon that can be spread to the whole plane by two vectors $\vec{k_1}=\left(\frac{a}{2},\frac{\sqrt{3}}{2}a\right)$ and $\vec{k_2}=(a,0)$ as shown in Figure \ref{fig4.3}.
  Inside the domain, there are six functions ($y_1,\cdots,y_6$)
 to be determined. Assume that each equilateral triangle has edgelength $a$. From Figure \ref{fig4.3}, we have first
  \begin{equation}
  -y_i''+q_i y_i=\la y_i\qquad (i=1,\ldots,6).\label{eq3.1}
  \end{equation}
  With the Neumann vertex conditions, we find that
 at $\bv_1$,
 \begin{equation}
 \left\{
 \begin{array}{l}
 y_1(0)=y_2(0)=y_3(0)=y_4(a)=y_5(a)=y_6(a)\\
y'_1(0)+y'_2(0)+y'_3(0)-y'_4(a)-y'_5(a)-y'_6(a)=0
 \end{array}
 \right.
 \label{eq3.5}
 \end{equation}
 And at vertices $\bv_2$, $\bv_3$, and $\bv_5$, by the Floquet-Bloch conditions,
 \begin{equation}
 \left\{
 \begin{array}{l}
 y_1(a/2)=e^{i\theta_1}y_4(a/2);\ y_1'(a/2)+e^{i\theta_1}y_4'(a/2)=0,\\
 y_2(a/2)=e^{i\theta_2}y_5(a/2);\ y_2'(a/2)+e^{i\theta_2}y_5'(a/2)=0,\\
 y_3(a/2)=e^{i(\theta_2-\theta_1)}y_6(a/2);\ y_3'(a/2)+e^{i(\theta_2-\theta_1)}y_6'(a/2)=0,
 \end{array}
 \right.
 \label{eq3.6}
 \end{equation}
  {This is because by periodicity, $q_1=q_4$, $q_2=q_5$ and $q_3=q_6$, and the edge $e_3$ is just $e_6$
 translated by the vector $\vec{k}_3=\vec{k}_2-\vec{k}_1$}. Hence
 Hence, by the uniqueness of solutions, \eqref{eq3.6}
 implies that
 $$
 y_4(x)=\al y_1(x);\quad y_5(x)=\be y_2(x);\quad y_6(x)=\al\be^{-1} y_3(x).
 $$
 where $\alpha=e^{-i\theta_1}$ and
 $\beta=e^{-i\theta_2}$.
 Since $y_j=A_jC_j+B_jS_j$, the equations in (\ref{eq3.5}) and (\ref{eq3.6}) implies
 \begin{equation}
 \left\{
 \begin{array}{l}
 A_1=A_2=A_3=\al(A_1 C_1+B_1 S_1)=\be(A_2C_2+B_2S_2)=\al^{-1}\be(A_3C_3+B_3S_3),\\
 -B_1-B_2-B_3+\al(A_1C_1'+B_1S_1')+\be(A_2C_2'+B_2S_2')+\al^{-1}\be(A_3C_3'+B_3S_3')=0.
 \end{array}
 \right. \label{eq3.65}
 \end{equation}
 This is a system of linear equations in $(A_1,A_2,A_3,B_1,B_2,B_3)$. There is a solution if and only if the following determinant  $\Phi(\rho)$ vanishes.
 \begin{center}
$\Phi(\rho) =
\begin{vmatrix}
 -1 & 1 & 0 & 0 & 0 & 0 &\\
 -1 & 0 & 1 & 0 & 0 & 0 &\\
 -1+\al C_1 & 0 & 0 & \al S_1 & 0 & 0 &\\
 -1 & \be C_2 & 0 & 0 & \be S_2 & 0 & \\
 -1 & 0 & \al^{-1}\be C_3 & 0 & 0 & \al^{-1}\be S_3 & \\
 \al C_1' & \be C_2' & \al^{-1}\be C_3' & -1+\al S_1' & -1+\be S_2' & -1+\al^{-1}\be S_3' &
\end{vmatrix}
$
\end{center}
 We expand the above determinant with the help of the symbolic software Mathematica. Then we group the terms in order of $\al^i\be^j$, and invoke the Lagrange identity $C_jS_j'-S_j C_j '=1$, to obtain
 \begin{eqnarray*}
 \Phi(\rho)&=& S_1 S_3(\be +\be^3)+S_2S_3 (\al^{-1}\be^2+\al\be^2)+S_1S_2(\al\be+\al^{-1}\be^3)-\be^2(S_1S_3(C_2+S_2')\\
  &&+S_2S_3(C_1+S_1')+S_1S_2(C_3+S_3')\\
  &=& -\rme^{-2i\theta_2} \left\{ (S_1S_2S_3)'+(C S S)_{1 2 3}-2S_1S_2\cos(\theta_2-\theta_1)-2S_2S_3\cos\theta_1\right.\\
  && \quad\left.-2S_1S_3\cos\theta_2\right\}.
 \end{eqnarray*}
 In the above, we used some simple identities $\be^2=\rme^{-2i\theta_2}$, and $\al+\al^{-1}=2\cos\theta_1$, while $\be +\be^{-1}=2\cos\theta_2$.
The notation $(CSS)_{123}$ was defined in \eqref{eq1.3}.
 \newtheorem{th2.1}{Theorem}[section]
  \begin{th2.1}
  \label{th2.1}
  For the periodic quantum graph associated with triangular tiling, the characteristic function is given by
  \begin{eqnarray*}
  \Phi(\rho)&=&-\rme^{-2i\theta_2} \left\{ (S_1S_2S_3)'+(C S S)_{1 2 3}-2S_1S_2\cos(\theta_2-\theta_1)-2S_2S_3\cos\theta_1\right.\\
  && \quad\left.-2S_1S_3\cos\theta_2\right\}.
  \end{eqnarray*}
  \end{th2.1}

  Assume that $q$ is even, then by Theorem \ref{th1.2}(d), $S'=C$. Therefore
  $$
   \Phi(\rho)=-2\rme^{-2i\theta_2} \left( (S_1S_2S_3)'-S_1S_2\cos(\theta_2-\theta_1)-S_2S_3\cos\theta_1-S_1S_3\cos\theta_2\right).
   $$
  And if all the potentials are equal, we let $C_j:=C$ and $S_j:=S$. Then the dispersion relation will be
  $$
  S^2\,\left( 3S'-\cos(\theta_2-\theta_1)-\cos\theta_1-\cos\theta_2\right)=0.
  $$
  Consequently by \eqref{eq1.2}, the dispersion relation becomes
  $$
  S^2\,\left( 3 S'+1-4\cos(\frac{\theta_1}{2})\cos(\frac{\theta_2}{2})\cos(\frac{\theta_2-\theta_1}{2})\right)=0.
  $$
 \section{Elongated triangular tiling}
\begin{figure}[h]
\centering
\includegraphics[width=7cm]{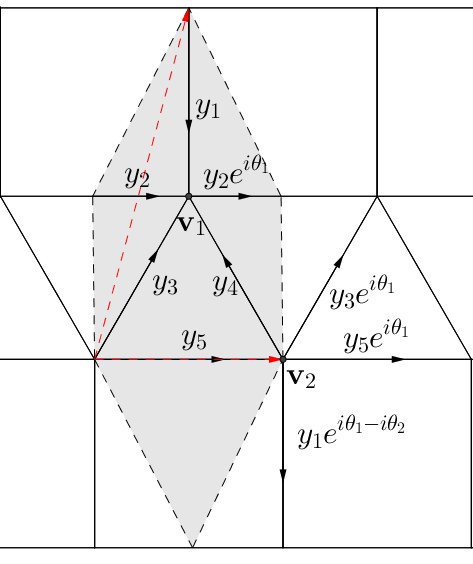}
\caption{Fundamental domain for elongated triangular tiling}
 \label{fig3.2}
\end{figure}

  \setcounter{equation}{0}\hskip0.25in
   For this periodic graph associated with elongated triangular
tiling with edgelength $a$, we choose a fundamental domain as a hexagon as shown in Fig.\ref{fig3.2} that can be spread to the whole plane
by two vectors $\vec{k}_1= {(a,0)}$ and $\vec{k}_2=\left(\frac{a}{2},(1+\frac{\sqrt{3}}{2})a\right)$.  {Thus $\displaystyle \vec{k}_3=\vec{k}_1-\vec{k}_2=(\frac{a}{2},-(1+\frac{\sqrt{3}}{2})a)$}.
 Realizing the Neumann vertex condition at $\bv_1$,
\begin{equation}  \label{1}
\begin{cases}
& y_1(a)=y_2(a)=y_2(0)e^{i\theta_1}=y_3(a)=y_4(a) \\
& y'_1(a)+y'_2(a)-y'_2(0)e^{i\theta_1}+y'_3(a)+y'_4(a)=0
\end{cases}
.
\end{equation}
Also, at the vertex $\bv_2$, we have
\begin{equation}  \label{2}
\begin{cases}
& y_4(0)=y_5(a)=y_3(0)e^{i\theta_1}=y_5(0)e^{i\theta_1}=y_1(0)e^{i\theta_1-i\theta_2}
\\
& y'_4(0)-y'_5(a)+y'_3(0)e^{i\theta_1}+y'_5(0)e^{i\theta_1}+y'_1(0)e^{i\theta_1-i\theta_2}=0%
\end{cases}%
.
\end{equation}
Since $y_j=A_jC_j+B_jS_j$, we can rewrite $(\ref{1})$ and $(\ref{2})$ as a
linear system of $(A_1,\cdots, A_5)$ and $(B_1,\cdots, B_5)$ with $
\widetilde{\alpha}=e^{i\theta_1}$ and $\beta=e^{-i\theta_2}$. Thus,
 $$
\begin{cases}
& A_1C_1+B_1S_1=A_2C_2+B_2S_2=A_2\widetilde{\alpha}=A_3C_3+B_3S_3=A_4C_4+B_4S_4 \\
& (A_1C'_1+B_1S'_1)+(A_2C'_2+B_2S'_2)-B_2\widetilde{\alpha}+(A_3C'_3+B_3S'_3)+(A_4C'_4+B_4S'_4)=0 \\
& A_4=A_5C_5+B_5S_5=A_3\widetilde{\alpha}=A_5\widetilde{\alpha}=A_1\widetilde{\alpha}\beta \\
& B_4-(A_5C'_5+B_5S'_5)+B_3\widetilde{\alpha}+B_5\widetilde{\alpha}+B_1\widetilde{\alpha}\beta=0
\end{cases}
.
 $$
So the characteristic equation $\Phi(\rho)$ is a determinant of $10\times 10$
matrix\newline

\begin{center}
$\Phi(\rho) =
\begin{vmatrix}
C_1 & -\widetilde{\alpha} & 0 & 0 & 0 & S_1 & 0 & 0 & 0 & 0 &  \\
0 & C_2-\widetilde{\alpha} & 0 & 0 & 0 & 0 & S_2 & 0 & 0 & 0 &  \\
0 & -\widetilde{\alpha} & C_3 & 0 & 0 & 0 & 0 & S_3 & 0 & 0 &  \\
0 & -\widetilde{\alpha} & 0 & C_4 & 0 & 0 & 0 & 0 & S_4 & 0 &  \\
C'_1 & C'_2 & C'_3 & C'_4 & 0 &S'_1 & S'_2-\widetilde{\alpha} & S'_3 & S'_4 & 0 &\\
0 & 0 & 0 & -1 & C_5 & 0 & 0 & 0 & 0 & S_5 &  \\
0 & 0 & \widetilde{\alpha} & -1 & 0 & 0 & 0 & 0 & 0 & 0 &  \\
0 & 0 & 0 & -1 & \widetilde{\alpha} & 0 & 0 & 0 & 0 & 0 &  \\
\widetilde{\alpha}\beta & 0 & 0 & -1 & 0 & 0 & 0 & 0 & 0 & 0 &  \\
0 & 0 & 0 & 0 & -C'_5 & \widetilde{\alpha}\beta & 0 & \widetilde{\alpha} & 1 &
\widetilde{\alpha}-S'_5
\end{vmatrix}%
$
\end{center}

The above matrix is complicated to handle by using manual computation.
However, with the help of Mathematica, and the Lagrange identity,
$C_{j}S_{j}'-S_{j}C_{j}'=1$, it still comes up with a long expression.
Grouping in terms of different powers of $\widetilde{\alpha}$ and $\widetilde{\beta}$, we obtain
 $$
 \Phi (\rho )=\Psi_1+\Psi_2+\Psi_3+\Psi_4,
 $$
 where, as $\widetilde{\alpha}^2+1=2\widetilde{\alpha}\, \cos\theta_1$, and $1+\widetilde{\alpha}^2\beta^2=2\widetilde{\alpha}\beta\,\cos(\theta_1-\theta_2)$,
  $$
 \Psi_1 = -(\widetilde{\al}^6+\widetilde{\al}^2)\be S_1S_3S_4-2\widetilde{\alpha}^4\beta S_1S_3S_4
  = -4\widetilde{\alpha}^4\beta\cos^2 \th_1\, S_1S_3S_4,
  $$
 while
 \begin{eqnarray*}
 \Psi_2 &=& (\widetilde{\alpha}^5+\widetilde{\alpha}^3)\beta\left( S_1'S_2S_3S_4+S_2'S_1S_3S_4+S_3'S_1S_2S_4+S_4'S_1S_2S_3+S_5'S_1S_3S_4\right. \\
    &&\left. +C_1S_3S_4S_5+C_2S_1S_3S_4+C_3S_1S_4S_5+C_4S_1S_3S_5+C_5S_1S_3S_4\right)\\
    &=& 2\al^4\be \cos\theta_1\, \left((S_1S_2S_3S_4)'+(C S S S)_{1345}+S_1S_3S_4(C_2+S_5')\right) \ ;\\
 \Psi_3 &=& (\widetilde{\alpha}^5+\widetilde{\alpha}^3)\beta S_1S_2S_5
+(\widetilde{\alpha}^3+\widetilde{\alpha} ^5\beta ^2) S_2S_3S_5+\widetilde{\alpha}^4(1+\be^2) S_2S_4S_5\\
  &=& 2\widetilde{\alpha}^4\be\, \left(S_1 S_2 S_5\cos\theta_1+S_2 S_3 S_5\cos(\theta_1-\theta_2)+S_2 S_4 S_5\cos\theta_2 \right).
    \end{eqnarray*}
 The last part is most tedious.  But we follow the pattern of the above main terms and succeed in
 grouping them  as follows:
 \begin{eqnarray*}
 \Psi_4 &=&
     -\widetilde{\alpha} ^4\beta\left[ S_1'S_5'S_2S_3S_4+S_2'S_5'S_1S_3S_4+S_3'S_5'S_1S_2S_4\right. \\
    &&\left. +S_4'S_5'S_1S_2S_3+C_1'S_2S_3S_4S_5+C_3'S_1S_2S_4S_5+C_4'S_1S_2S_3S_5+S_1'C_3S_2S_4S_5+S_1'C_4S_2S_3S_5\right. \\
  &&\left. +S_1'C_5S_2S_3S_4+S_2'C_1S_3S_4S_5+S_2'C_3S_1S_4S_5+S_2'C_4S_1S_3S_5+S_2'C_5S_1S_3S_4+S_3'C_1S_2S_4S_5\right. \\
   &&\left. +S_3'C_4S_1S_2S_5+S_3'C_5S_1S_2S_4+S_4'C_1S_2S_3S_5+S_4'C_3S_1S_2S_5+S_4'C_5S_1S_2S_3+S_5'C_2S_1S_3S_4\right.\\
   &&\left. +C_1C_2S_3S_4S_5+C_2  C_3S_1S_4S_5+C_2C_4S_1S_3S_5+C_2C_5S_1S_3S_4\right] \\
   &=& {-\widetilde{\alpha} ^4\beta}\left[(S_1S_2S_3S_4)'(C_5+S_5')+(C S S S)_{1345}(S_2'+C_2)+S_1S_2S_5((C_4S_3)'+(C_3 S_4)')\right.\\
  &&\left.+S_2S_3S_5(C_4S_1'+(C_1S_4)')+S_2S_4S_5( {S'_1C_3}+C_1S_3') +S_1S_3S_4(C_2S_5'-S_2'C_5)\right]
  \end{eqnarray*}
Therefore, we have the following formula for the characteristic function:
 \newtheorem{th3.1}{Theorem}[section]
 \begin{th3.1}
 \label{th3.1}
For the periodic quantum graph of elongated triangular tiling with length $a$, we have
\begin{eqnarray*}
\Phi(\rho)&=&-\widetilde{\alpha}^4\beta\{ (S_1S_2S_3S_4)'(C_5+S_5'-2\cos\theta_1)+(C S S S)_{1345}(S_2'+C_2-2\cos\theta_1) \\
   &&+S_1S_2S_5((C_4S_3)'+(C_3 S_4)'-2\cos\theta_1)+S_2 S_3 S_5((C_1S_4)'+S_1'C_4-2\cos(\theta_1-\theta_2))\\
 &&+S_2S_4S_5( {S'_1C_3}+C_1S_3'-2\cos\theta_2)+S_1S_3S_4(C_2S_5'-S_2'C_5-2(C_2+S_5')\cos\theta_1 {+}4\cos^2 \theta_1)\}
\end{eqnarray*}
\end{th3.1}

  Note that in Theorem \ref{th3.1}, the involved quantum graph has different potential functions on distinct edges. If all the potentials are equal,
  then we let $C_j:=C$ and $S_j:=S$. The dispersion relation is given by
 \begin{eqnarray*}
  0&=&S^3\{2\cos\theta_1+2\cos\theta_2+2\cos(\theta_1-\theta_2)+3+10S'\cos\theta_1+10C\cos\theta_1-4\cos^2\theta_1\\
 &&-17S'C-4(S')^2-4C^2\}
\end{eqnarray*}
If we further assume that the potentials are even, then by Theorem~\ref{th1.2}(d),
 the dispersion relation becomes
 \begin{eqnarray*}
 0 &=&S^3\{25(S')^2-20\cos\theta_1S'-2\cos\theta_1-2\cos\theta_2-2\cos(\theta_1 - \theta_2)+4\cos^2\theta_1-3\}\\
 &=& S^3\{25(S')^2-20\cos\theta_1S'-8\cos(\frac{\theta_1}{2})\cos(\frac{\theta_2}{2})\cos(\frac{\theta_1-\theta_2}{2})+4\cos^2\theta_1-1\},
 \end{eqnarray*}
  by \eqref{eq1.2}. This proves Theorem \ref{th1.0}(b).
 \\[0.2in]
 \section{Truncated square tiling}
 \setcounter{equation}{0}\hskip0.25in
  For this periodic quantum graph associated with truncated square
 tiling with edgelength $a$, we choose a fundamental domain as a square
 spanned by two vectors $\vec{k}_1=\left(0,(1+\sqrt{2})a\right)$ and $\vec{k}_2=\left( (1+\sqrt{2})a,0\right)$ as shown in the figure below.
\begin{figure}[h!]
 \centering\includegraphics[width=6cm,height=6cm]{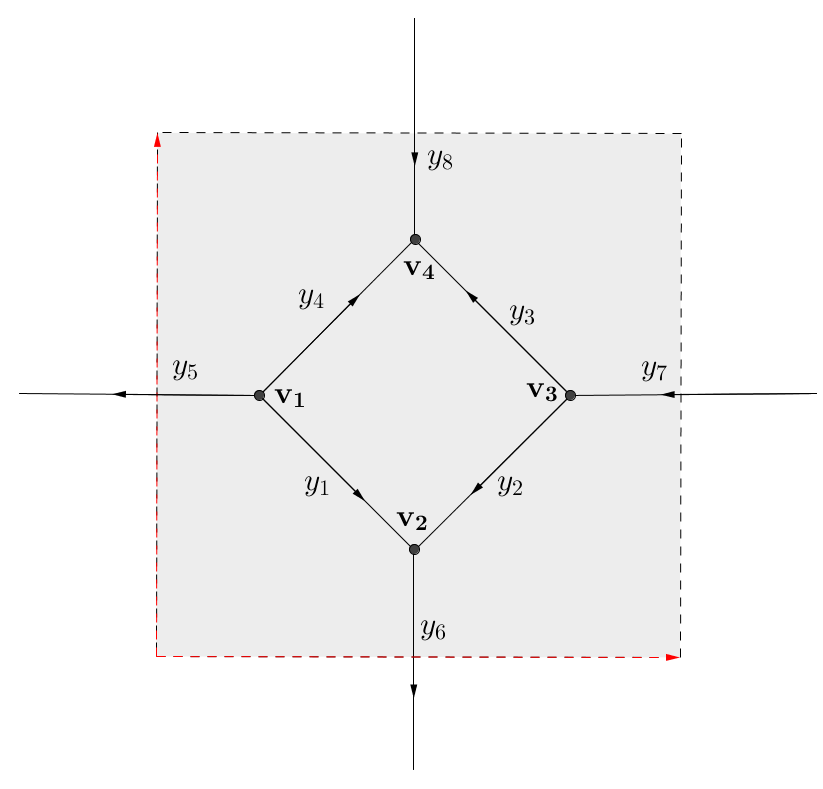}
 \caption{Fundamental domain for truncated square tiling}
 \label{fig4.1}
 \end{figure}
We first note that by the Floquet-Bloch conditions, $y_7(x)=\rme^{i\theta_1}y_5(x)$, and $y_8(x)=\rme^{i\theta_2}y_6(x)$.
 Next  the Neumann vertex conditions at vertex $\bv_1$ imply
 $$
\begin{cases}
& y_1(0)=y_4(0)=y_5(0) \\
& y'_1(0)+y'_4(0) +y'_5(0)=0
\end{cases}
.
  $$
 At the vertex \textbf{v}$_2$,
 $$
\begin{cases}
& y_1(a)=y_2(a)= y_6(0)
\\
& y'_1(a)+y'_2(a)-y'_6(0)=0%
\end{cases}
.
  $$
 At vertex \textbf{v}$_3$, they become
 $$
\begin{cases}
& y_2(0)=y_3(0)= \rme^{i\theta_1}y_5(a)
\\
& y'_2(0)+y'_3(0)-\rme^{i\theta_1}y_5(a)=0%
\end{cases}
 ;
  $$
  and
  $$
\begin{cases}
& y_3(a)=y_4(a)= \rme^{i\theta_2}y_6(a)
\\
& y'_3(a)+y'_4(a)+\rme^{i\theta_2}y'_6(a)=0
\end{cases}
.
  $$
Since $y_j=A_jC_j+B_jS_j$, we can rewrite the above equations as a
linear system of $(A_1,\cdots, A_6)$ and $(B_1,\cdots, B_6)$ with $%
\widetilde{\alpha}=e^{i\theta_1}$ and $\widetilde{\beta}=e^{i\theta_2}$. Thus,
 $$
\begin{cases}
 &A_1=A_4=A_5\\
 &B_1+B_4+B_5=0\\
 &A_6=A_1 C_1+B_1S_1=A_2C_2+B_2S_2\\
 &B_6-A_1C_1'-B_1 S_1'-A_2C_2'-B_2 S_2=0\\
 &A_2=A_3=\rme^{i\theta_1}(A_5 C_5+B_5C_5)\\
 &B_2+B_3-\rme^{i\theta_1}(A_5 C_5'+B_5 S_5')=0\\
 & \rme^{i\theta_2}(A_6 C_6+B_6 S_6)=A_3 C_3+B_3 S_3=A_4 C_4+B_4 S_4\\
 & \rme^{i\theta_2}(A_6 C_6'+B_6 S_6') + A_3 C_3'+B_3 S_3' +A_4 C_4'+B_4 S_4'=0
\end{cases}
.
 $$
So the characteristic equation $\Phi(\rho)$ is a determinant of $12\times 12$
matrix\newline

\begin{center}
$\Phi(\rho) =
\scriptsize{\begin{vmatrix}
1 & 0 & 0 & -1 & 0 &0 & 0 & S_1 & 0 & 0 & 0 & 0 &  \\
1 & 0 & 0  & 0 & -1 & 0 & 0 & 0 & 0 & 0 & 0 & 0 &  \\
0 & 0 & 0 & 0 & 0 & 0 & 1 & 0 & 0 & 1 & 1 & 0 & \\
C_1 & 0 & 0 & 0 & 0 & -1 & S_1 & 0 & 0 & 0 & 0 & 0 &\\
0 & C_2 & 0 & 0 & 0 & -1 & 0 & S_2 & 0 & 0 & 0 & 0\\
C'_1 & C'_2 & 0 & 0 & 0  & 0 &S'_1 & S'_2 & 0 & 0 & 0 & -1 &\\
0 & 1 &  -1 & 0 & 0 & 0 & 0 & 0 & 0 & 0 & 0 & 0 &  \\
0 & -1 & 0 & 0 & \al C_5 & 0 & 0 & 0 & 0 & 0 & \al S_5 & 0 &\\
0 & 0 & 0 & 0 & \al C_5' & 0 & -1 & -1 & 0 & 0 & \al S_5' & 0 &\\
0 & 0 & -C_3 & C_4 & 0 & 0 & 0 & 0 &-S_3 & S_4 & 0 & 0 &  \\
0 & 0 & -C_3 & 0 & 0 & \be C_6 & 0 & 0 & -S_3 & 0 & 0 & \be S_6 &\\
0 & 0 & C_3' & C_4' & 0 & \be C_6' & 0 & 0 & S_3' & S_4' & 0 & \be S_6' &
\end{vmatrix}}
$
\end{center}

The above matrix is complicated to handle by using manual computation.
However, with the help of Mathematica, we manage to simplify the expression using the Lagrange identity,
$C_{j}S_{j}'-S_{j}C_{j}'=1,\forall j=1,\cdots ,5$.  It still comes up with a long expression.
Grouping in terms of different powers of $\widetilde{\alpha}$ and $\widetilde{\beta}$, we obtain
 $$
 -\Phi (\rho )=\Psi_1+\Psi_2+\Psi_3+\Psi_4,
 $$
 where
 $$
 \Psi_1=(1+\widetilde{\alpha}^2\widetilde{\beta}^2) S_1 S_3+(\widetilde{\alpha}^2+\widetilde{\beta}^2)S_2S_4
 2 S_2S_4\cos(\theta_1-\theta_2)+2 S_1S_3\cos(\theta_1+\theta_2),
 $$
 while
 \begin{eqnarray*}
 \Psi_2&=& (1+\widetilde{\alpha}^2)\widetilde{\beta} (C_6S_1S_2+S_1S_6S_2'+S_2S_6S_1'+S_4S_6S_3' +S_3S_6S_4'+S_3S_4S_6')\\
  &=& 2\widetilde{\alpha}\widetilde{\beta} \cos\theta_1\, \left((S_3S_4S_6)'+S_6(S_1S_2)'+C_6S_1S_2\right)\,\\
  \Psi_3 &=&\widetilde{\alpha}(1+\widetilde{\beta}^2)(C_5S_1S_4 + C_4S_1S_5+C_3S_2S_5+C_2S_3S_5+C_1S_4S_5+S_2S_3S_5')\\
  &=&  2\widetilde{\alpha}\widetilde{\beta}\cos \theta_2\, \left( (CSS)_{145}-S_5(CS)_{23}-S_2S_3S_5'\right).
  \end{eqnarray*}
  Finally, $\Psi_4$ is the most tedious part.
\begin{footnotesize}
\begin{eqnarray*}
 \Psi_4&=& \widetilde{\alpha}\widetilde{\beta}(-2 S_5S_6
 +C_1C_6S_2S_4S_5C_3'+S_2 S_4S_5S_6C_1'C_3'+C_3C_6S_1S_2S_5C_4'+C_2C_6S_1S_3S_5C_4'+S_1S_3S_5S_6C_2'C_4'\\
 && +C_1C_3S_2S_4S_5C_6'C_3S_2S_5S_6C_4'S_1'+C_2S_3S_5S_6C_4'S_1'+C_1S_4S_5S_6C_3'S_2'+C_3S_1S_5S_6C_4'S_2'+C_1C_2C_6S_4S_5S_3'\\
 && +C_2S_4S_5S_6C_1'S_3'
 +C_1S_4S_5S_6C_2'S_3'+C_6S_1S_2S_4C_5'S_3'+S_2S_4S_6C_5'S_1'S_3'+S_1S_4S_6C_5'S_2'S_3'+C_1C_3C_6S_2S_5S_4'\\
  && +C_1C_2CC_6S_3S_5S_4'+C_3S_2S_5S_6C_1'S_4'+C_2S_3S_5S_6C_1'S_4'
 +C_1S_3S_5S_6C_2'S_4'+C_6S_1S_2S_3C_5'S_4'+S_2S_3S_6C_5'S_1'S_4'\\
 && +C_1C_3S_5S_6S_2'S_4'+S_1S_3C_6C_5'S_2'S_4'+C_6S_1S_2S_3C_4'S_5'+C_1S_2S_3S_4C_6'S_5'+S_2S_3S_6C_4'S_1'S_4'+S_1S_3S_6C_4'S_2'S_5'\\
 && +C_1C_6S_2S_4S_3'S_5'+S_2S_4S_6C_1'S_3'S_5'+C_1S_4S_6S_2'S_3'S_5'C_1C_6S_2S_3S_4'S_5'+S_2S_3S_6C_1'S_4'S_5'+C_1S_3S_6S_2'S_4'S_5'\\
 && +C_3S_2S_4S_5C_1'S_6'+C_2S_3S_4S_5C_1'S_6'+C_1S_3S_4S_5C_2'S_6'
 +S_2S_3S_4C_5'S_1'S_6'+C_1C_3S_4S_5S_2'S_6'+S_1S_3S_4C_5'S_2'S_6'\\
 && +S_2S_3S_4C_1'S_5'S_6'+C_1S_3S_4S_2'S_5'S_6'+C_4\, (C_2S_1S_3S_5C_6'+S_2S_5S_6C_3'S_1'+S_1S_5S_6C_3'S_2'+S_1S_5S_6C_2'S_3'\\
 && +C_2S_5S_6S_1'S_3'+S_1S_2S_3C_6'S_5'+S_2S_6S_1'S_3'S_5'
 +S_1S_6S_2'S_3'S_5'+C_6S_1(C_2S_5S_3'+S_2(S_5C_3'+S_3'S_5'))\\
 && +S_1S_3S_5C_2'S_6'+C_2S_3S_5S_1'S_6'+S_2S_3S_1'S_5'S_6'+S_1S_3S_2'S_5'S_6'+C_3S_5(S_2S_1'S_6'+S_1(S_2C_6'+S_2'S_6')))\\
 && + C_5\, (C_2S_1S_3S_4C_6'+S_2S_4S_6C_3'S_1'+S_1S_4S_6C_3'S_2'+S_1S_4S_6C_2'S_3'+C_2S_4S_6S_1'S_3'+ S_1S_3S_6C_2'S_4'\\
 && +C_2S_3S_6S_1'S_4'+C_6S_1(S_2(S_4C_3'+C_3S_4')+C_2(S_4S_3'+S_3S_4'))+S_1S_3S_4C_2'S_6'+C_2S_3S_4S_1'S_6'\\
 && +C_3(S_2S_1'(S_6S_4'+S_4S_6')+S_1(S_2S_4C_6'+S_2'(S_6S_4'+S_4S_6')))))
 \end{eqnarray*}
 \end{footnotesize}
 We follow the pattern of the above main terms and succeed in
 grouping them, as follows:
 \begin{footnotesize}
\begin{eqnarray*}
 \Psi_4
 &=&(S_5S_2S_3+S_5(CS)_{23})\, (C_1'(S_4S_6)'+S_1'(C_4S_6)'+S_1(C_4C_6)'+C_1(S_4C_6)')+(CSS)_{145}\, (S_2'(C_3S_6)'+C_2'(S_3S_6)'\\
  &&+C_6(CS)_{2 3'})+(C_6S_1S_2+S_6(S_1S_2)')\, (C_5'(S_3S_4)'+C_3S_4'C_5+S_3'C_4S_5') + (S_3S_4S_6)'\, (C_1S_2'S_5'+S_1'C_2C_5)\\
  &&+S_5S_6\, ((SC)_{14'}(SC)_{3 2'}+(SC)_{4 1'}(SC)_{2 3'})+S_2S_4\, ( (S_1'S_6'+S_1C_6')(C_3C_5+S_3C_5')\\
  &&+S_6(C_1'S_3'S_5'+S_1'C_3'C_5)+C_1C_6S_3'S_5')
  +S_1S_3\, (S_2'S_6'(CS)_{4 5'}+C_2C_5(S_4S_6)'+S_6(C_2'S_4'C_5+S_2'C_4'S_5').
  \end{eqnarray*}
 \end{footnotesize}
 \hskip0.25in
Therefore, we have the following formula for the characteristic function:
 \newtheorem{th4.1}{Theorem}[section]
 \begin{th4.1}
 \label{th4.1}
For the periodic quantum graph of truncated square tiling with length $a$, we have
\begin{eqnarray*}
\Phi(\rho)&=&\rme^{i\theta_1+i\theta_2}\, \left\{ (S_5'S_2S_3+S_5(CS)_{23})\, (C_1'(S_4S_6)'+S_1'(C_4S_6)'+S_1(C_4C_6)'+C_1(S_4C_6)'-2\cos\theta_2)\right.\\
 && +(CSS)_{145}\, (S_2'(C_3S_6)'+C_2'(S_3S_6)'+C_6(CS)_{2 3'}-2\cos\theta_2)\\
 && +(C_6S_1S_2+S_6(S_1S_2)')\, (C_5'(S_3S_4)'+C_3S_4'C_5+S_3'C_4S_5'-2\cos\theta_1)+(S_3S_4S_6)'\, (C_1S_2'S_5'\\
 && +S'_1C_2C_5-2\cos \theta_1)+ S_5S_6\, ((SC)_{14'}(SC)_{3 2'}+(SC)_{4 1'}(SC)_{2 3'}-2)\\
  && +S_2S_4\, ( (S_1'S_6'+S_1C_6')(C_3C_5+S_3C_5')+S_6(C_1'S_3'S_5'+S_1'C_3'C_5)+C_1C_6S_3'S_5'-2\cos(\theta_1-\theta_2))\\
  && \left. +S_1S_3\, (S_2'S_6'(CS)_{4 5'}+C_2C_5(S_4C_6)'+S_6(C_2'S_4'C_5+S_2'C_4'S_5')-2\cos(\theta_1+\theta_2))\right\}.
\end{eqnarray*}
\end{th4.1}
\noindent
  Note that in Theorem~\ref{th4.1}, the involved quantum graph has different potential functions on distinct edges. If all the potentials are equal, then the dispersion relation becomes
 \begin{eqnarray*}
  0&=& S^2\left\{ C S'(18 S'^2+45 C S'+18 C^2)-11 S'^2-32 C S-11 C^2+1-2(5S'+C)\cos \theta_1\right.\\
  && \left. -2(S'+5C)\cos\theta_2-2\cos(\theta_1-\theta_2)-\cos(\theta_1+\theta_2)\right\}
\end{eqnarray*}
 If we further assume that the potential is even, then we may invoke
  a trigonometric identity to arrive at the following dispersion relation.
 $$
 S^2\left\{ 81 S'^4-54 S'^2-12S'(\cos\theta_1+\cos\theta_2)+1-4\cos\theta_1\, \cos\theta_2\right\} = 0.
 $$
 \vskip0.2in
  \section{Trihexagonal tiling}
  \setcounter{equation}{0}\hskip0.25in
  For this periodic graph associated with trihexagonal tiling with edgelength $a$, we choose a fundamental domain as a parallelogram spanned by two vectors
  $\vec{k}_1=(\sqrt{3} a, a)$, and $\vec{k}_2=(0,2a)$.
  \begin{figure}[h!]
 \centering\includegraphics[width=7cm,height=7cm]{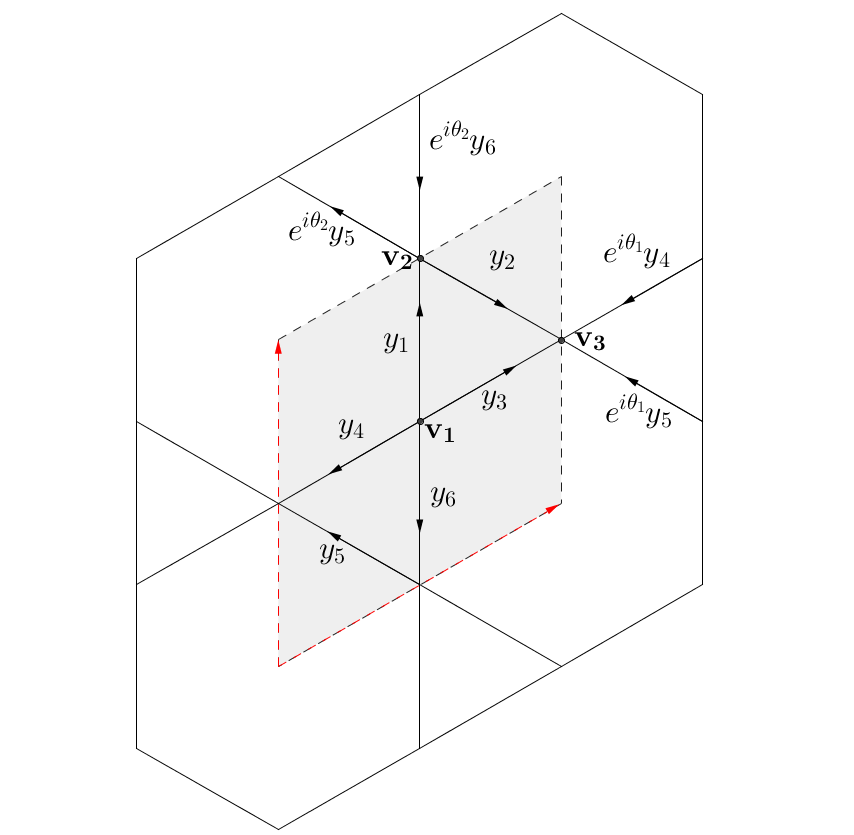}
 \caption{Fundamental domain for trihexagonal tiling}
 \label{fig5.1}
 \end{figure}
 Letting $y_i(x)=A_i C_i(x)+B_i S_i(x)$, $i=1,\ldots,6$. Using the Neumann vertex conditions, we have
  \begin{equation*}
  \left\{
  \begin{array}{l}
  A_1=A_3=A_4=A_6\\
  B_1+B_3+B_4+B_6=0\\
  A_2C_2+B_2S_2=A_3 C_3+B_3S_3=\rme^{i \theta_1} (A_4 C_4+B_4S_4)=\rme^{i \theta_1} (A_5 C_5+B_5S_5)\\
  A_2 C_2'+B_2S_2'+A_3 C_3'+B_3S_3'+ \rme^{i \theta_1}(A_4 C_4'+B_4S_4')
  + \rme^{i \theta_1}(A_5 C_5'+B_5S_5')=0\\
  A_2=\rme^{i \theta_2}A_5=\rme^{i \theta_2}(A_6 C_6+B_6S_6)=A_1 C_1+B_1S_1\\
 -B_2-\rme^{i \theta_2}B_5+\rme^{i \theta_2}(A_6 C_6'+B_6S_6')+A_1 C_1'+B_1S_1'=0
 \end{array}
 \right.
 \end{equation*}

   Let $\displaystyle  \widetilde{\alpha}=\rme^{i\theta_1}$, $\widetilde{\beta}=\rme^{i\theta_2}$. The characteristic equation $\Phi(\rho)$ is a determinant of another $12\times12$ matrix
 \begin{footnotesize}
 \begin{eqnarray*}
 \Phi(\rho)&=&
 \left|
 \begin{array}{cccccccccccc}
 1 & 0 & -1 & 0 & 0 & 0 & 0 & 0 & 0 & 0 & 0 & 0\\
 1 & 0 & 0 & -1 & 0 & 0 & 0 & 0 & 0 & 0 & 0 & 0\\
 1 & 0 & 0 & 0 & 0 & -1 & 0 & 0 & 0 & 0 & 0 & 0\\
 0 & 0 & 0 & 0 & 0 & 0 & 1 & 0 & 1 & 1 & 0 & 1\\
 0 & C_2 & -C_3 & 0 & 0 & 0 & 0 & S_2 & -S_3 & 0 & 0 & 0\\
 0 & C_2 & 0 & -\widetilde{\alpha} C_4 & 0 & 0 & 0 &S_2 & 0 & -\widetilde{\alpha} S_4 & 0 & 0\\
 0 & C_2 & 0 & 0&-\widetilde{\alpha} C_5 & 0 & 0 & S_2 & 0 & 0&-\widetilde{\alpha} S_5 & 0 \\
 0 & C_2'& C_3' & \widetilde{\alpha} C_4' & \widetilde{\alpha} C_5' & 0 &0 & S_2' & S_3' &\widetilde{\alpha} S_4' & \widetilde{\alpha} S_5' & 0 \\
 0 & 1 & 0 & 0 & -\widetilde{\beta} & 0 & 0 & 0 & 0 & 0 & 0 & 0\\
 0 & 1 & 0 & 0 & 0 &-\widetilde{\beta} C_6  & 0 & 0 & 0 & 0 & 0 & -\widetilde{\beta} S_6\\
 -C_1 & 1 & 0 & 0 & 0 & 0 & -S_1 & 0 & 0 & 0 & 0 & 0 \\
 C_1' & 0 & 0 & 0 & 0 & \widetilde{\beta} C_6' &  S_1' & -1 & 0 &0 & -\widetilde{\beta} & \widetilde{\beta} S_6'
 \end{array}
 \right|
 \end{eqnarray*}
 \end{footnotesize}

 Using the software Mathematica, we expand the determinant and group the terms in the order of $\widetilde{\alpha}^i\widetilde{\beta}^j$. After simplification
 \begin{eqnarray*}
 -\Phi(\rho) &=& (\widetilde{\alpha}\widetilde{\beta}^2+\widetilde{\alpha}^3\widetilde{\beta}^2)(S_1S_2S_4+S_3S_5S_6+S_1S_6(CS)_{25}+S_2S_5(S_1S_6)'+(\widetilde{\alpha}\widetilde{\beta}^3+\widetilde{\alpha}^3\widetilde{\beta})(S_1S_3S_5\\
 && +S_2S_4S_6+(CSSS)_{1346})+(\widetilde{\alpha}^2\widetilde{\beta}+\widetilde{\alpha}^2\widetilde{\beta}^3)(S_1S_4S_5+S_2S_3S_6+(S_2S_3S_4S_5)')-\widetilde{\alpha}^2\widetilde{\beta}^2\, \Delta.
 \end{eqnarray*}
 Here $\Delta$ is a very long algebraic expression consisting of around 70 terms.  Following the pattern in previous section, we obtain
 \begin{eqnarray*}
 \Delta &=&  -2(S_1S_2S_3+S_4S_5S_6)+(S_2 S_3 S_4 S_5)'(S_6C_1+S_1C_6)'+(C S S S)_{1346}(C_2S_5+C_5S_2)'\\
 && + (S_1S_6)'(S_2S_5)'(C_3S_4+C_4S_3)+(C_3S_4+C_4S_3)'\, (S_1S_6(C_2S_5+C_5S_2)+S_2S_5(S_1S_6)')\\
 &&+(S_3S_4)'(C_2S_5+C_5S_2)(C_1S_6)+C_6S_1).
 \end{eqnarray*}
 \newtheorem{th5.1}{Theorem}[section]
 \begin{th5.1}
 \label{th5.1}
 For the periodic quantum graph of trihexagonal tiling, the dispersion relation is given by
 \begin{eqnarray*}
 \Phi(\rho) &=& -\rme^{2i(\theta_1+\theta_2)}\left\{ 2((S_1 S_2 S_4+S_3 S_5 S_6)\cos\theta_1+(S_1 S_3 S_5+S_2 S_4 S_6) \cos(\theta_1-\theta_2)\right.\\
  &&+(S_1 S_4 S_5+S_2 S_3 S_6)\cos\theta_2+S_1 S_2 S_3+S_4 S_5 S_6)+(S_2 S_3 S_4 S_5)'\, (2\cos\theta_2-((C S)_{16})')\\
  && +(C S S S)_{1346}(2\cos(\theta_1-\theta_2)-((C S)_{25})') -(S_1S_6)'(S_2S_5)'(C S)_{34}\\
  && \left.-(S_3S_4)'(C S)_{25}( C S)_{16}+(S_1S_6 (C S)_{25}+S_2S_5(S_1S_6)')(2\cos\theta_1-((C S)_{34})')\right\}
  \end{eqnarray*}
  \end{th5.1}

  When all potentials are identical to $q$, then using Lagrange identity, it yields
  \begin{eqnarray*}
  \Phi(\rho)&=& -4\widetilde{\alpha}^2\widetilde{\beta}^2S^3 \left\{ 8C S'^2+8S'C^2)-C(3+2\cos(\theta_1-\theta_2)+\cos\theta_1)\right. \\
   &&\quad \left.-S'(3+2\cos\theta_2+\cos\theta_1)-\cos\theta_1-\cos\theta_2-\cos(\theta_1-\theta_2)-1\right\}
  \end{eqnarray*}
 And when $q$ is even, the dispersion relation becomes
 \begin{eqnarray*}
 0 &=&S^3\left( 4 S'^3-S'(1+2\, \cos(\frac{\theta_1}{2})\cos(\frac{\theta_1}{2})\cos(\frac{\theta_1-\theta_2}{2}))-\cos(\frac{\theta_1}{2})\cos(\frac{\theta_1}{2})\cos(\frac{\theta_1-\theta_2}{2})\right)\\
  &=&  {S^3(2S'+1)\, \left( 2S'^2-S'-\cos(\frac{\theta_1}{2})\cos(\frac{\theta_1}{2})\cos(\frac{\theta_1-\theta_2}{2})\right)}
 \end{eqnarray*}

\begin{figure}[h!]
\begin{minipage}{.5\linewidth}
\centering
\subfloat[]{\label{main:a}\includegraphics[width=4cm,height=4cm]{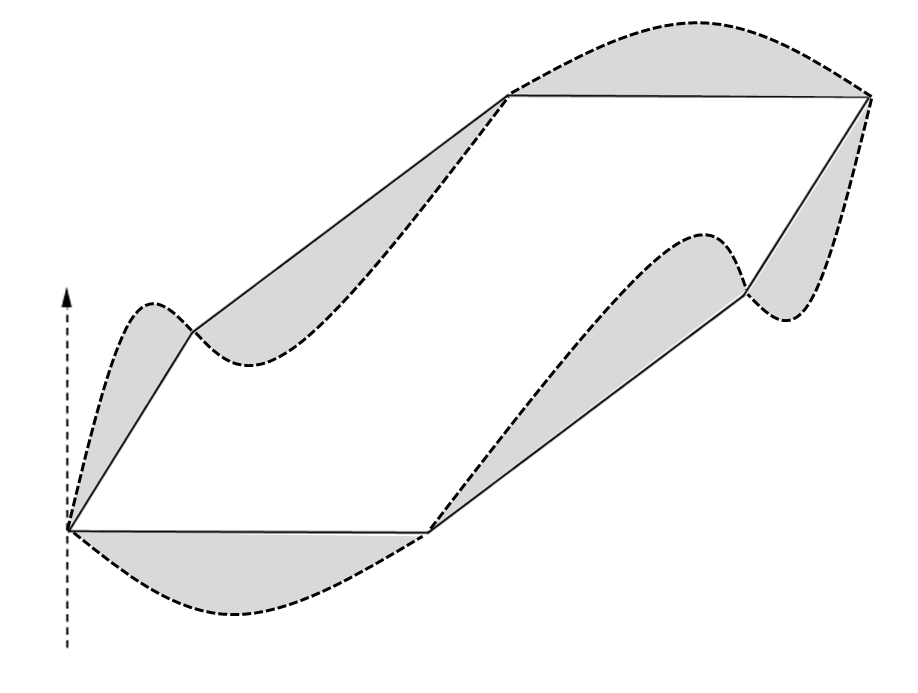}}
\end{minipage}%
\begin{minipage}{.5\linewidth}
\centering
\subfloat[]{\label{main:b}\includegraphics[width=4cm,height=4cm]{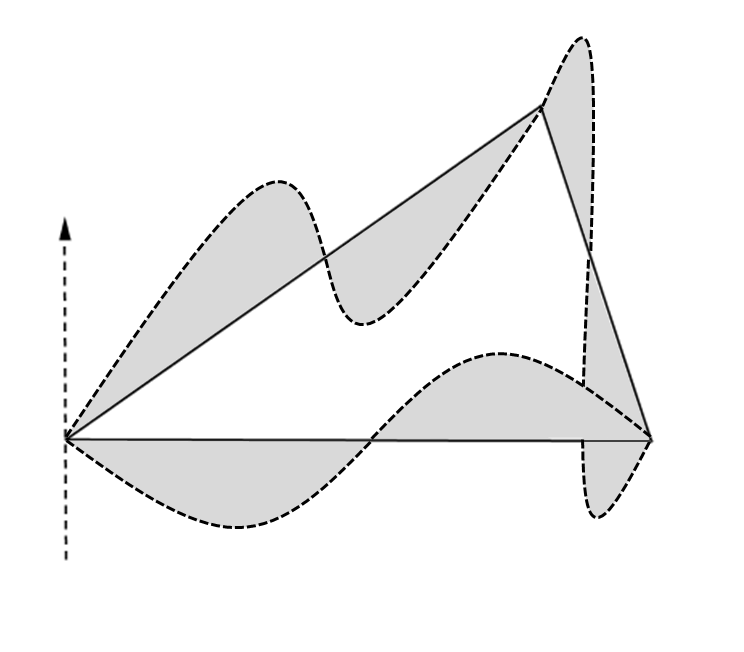}}
\end{minipage}\par\medskip
\centering
\subfloat[]{\label{main:c}\includegraphics[width=12cm,height=4cm]{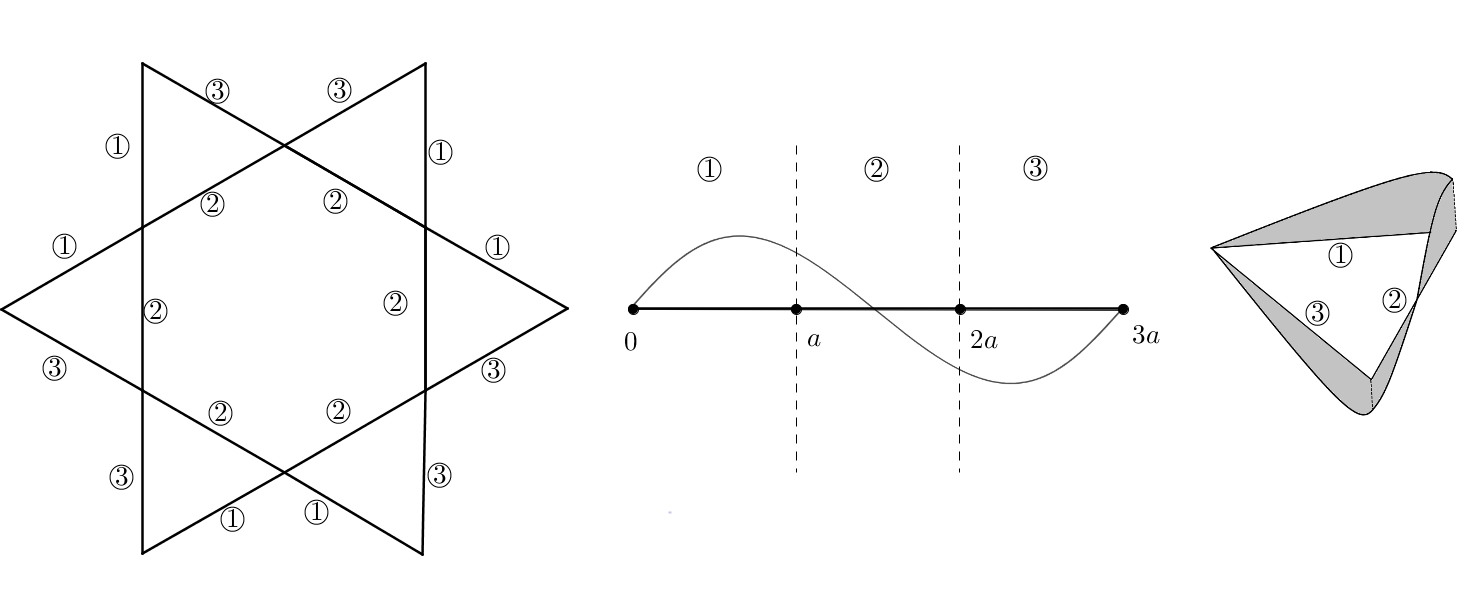}}
\caption{Eigenfunctions for trihexagonal tiling}
 \label{fig3.1}
 \end{figure}

  Similar as the case of hexagonal tiling \cite{KP,KL}, the equation $S(a,\rho)=0$ defines a point spectrum with infinite number of eigenfunctions. A typical example of these eigenfunctions is a repetition of the function $\vp_2$ for the Dirichlet-Dirichlet boundary conditions around any basic triangular unit, while vanishing elsewhere.  Another example is a repetition of $\vp_1$ around any basic hexagonal unit.
 There are infinite number of such eigenfunctions on the periodic quantum graph.  {Besides the two examples as} the following shaded functions on a hexagon  {(type (a)} or a triangle  {(type (b)}, extended by zero to the rest of the graph,  {another class of eigenfunctions (type (c)) are also found} (see Fig.\ \ref{fig3.1}).\\[0.1in]
 \hskip0.25in
  {Let us explain the class of type (c) eigenfunctions. By the theory of Hill's discriminant \cite[Theorem 2.3.1]{E73}, the condition $S'(a,\rho)=-1/2$ defines a sequence of eigenpairs $(\la_n, u_n)$ on the interval $[0,a]$.
 This sequence of solutions $u_n$, if extended, has $n$ periods in $(0,3a)$. The graphs, on each subintervals are different, and are names as graphs $\textcircled{1}, \textcircled{2}, \textcircled{3}$.
 They represent Dirichlet-Dirchlet eigenfunctions on $[0,3a]$ because from Theorem~\ref{th1.2} (a),(b) and (d), with $\la_n=\rho_n^2$,
 \begin{eqnarray*}
 S'(2a,\rho_n)&=&C(2a,\rho_n)=2\, S'(a,\rho_n)^2-1=-\frac{1}{2};\\
  S(2a,\rho_n)&=&2\, S(a,\rho_n)S'(a,\rho_n)=-S(a,\rho_n).
 \end{eqnarray*}
 Notice that  if we fold the graph to base it on a triangle with sidelength $a$, and place six of them on the sides of a hexagon, allowing the function to vanish
  on the rest of the graph (see Fig.\ \ref{fig3.1}), we obtain a sequence of eigenfunctions for the quantum graph with eigenvalues $\la_n$'s. The multiplicity for each
  eigenvalue $\la_n$ is infinite.}
 \vskip0.2in
 \section{More analysis on the spectrum}
 \setcounter{equation}{0}\hskip0.25in
 In this section, we study the spectra of these periodic quantum graphs in more detail.  In short, we try to understand more about the point spectrum
 $\sig_p$ and the asymptotic behavior of the absolutely continuous spectrum $\sig_{ac}$. Define $\bfN_k=\{ 1,2,\ldots,k\}$, for any $k\in \bfN$.
 \newtheorem{th6.1}{Theorem}[section]
 \begin{th6.1}
 \label{th6.1}
 \begin{enumerate}
 \item[(a)] The point spectrum of the periodic quantum graph associated with triangular tiling satisfies
 $$
  \sig_p(H_T)= \bigcup_{i,j\in \bfN_3}\{ \rho^2\in\bfR:\ S_i(a,\rho)=S_j(a,\rho)=0\}
 $$
 \item[(b)] The point spectrum of the periodic quantum graph associated with elongated triangular tiling satisfies
 $$
 \sig_{p}(H_{eT})= \bigcup_{i,j,k\in \bfN_5}\{ \rho^2\in\bfR:\ S_i(a,\rho)=S_j(a,\rho)=S_k(a,\rho)=0\}
 $$
 \item[(c)] The point spectrum of the periodic quantum graph associated with truncated square tiling satisfies:\
  $\la=\rho^2$ lies in $\sig_p(H_{trS})$ if and only if $\rho$ satisfies
  \begin{eqnarray}
  2S_1 S_3&=&2 S_2S_4\ =\ C_6S_1S_2+S_6(S_1S_2)'+(S_3S_4S_6)'\nonumber\\
    &=& S_5'S_2S_3+S_5(CS)_{23}+(C S S)_{145}\ =\ 0 \ ,\label{eq6.1}
    \end{eqnarray}
    and
    \begin{eqnarray}
    0 &=&(C S S)_{145}\, \left( S_2'(C_3S_6)'+C_2'(S_3S_6)'+C_6(C S)_{2 3'}-C_1'(S_4S_6)'-S_1'(C_4S_6)'-S_1(C_4C_6)'\right.
    \nonumber\\
     &&\left. -C_1(S_4C_6)'\right)+(S_3 S_4 S_6)'\, \left( C_1 S_2'S_5'+S_1'C_2C_5-C_5'(S_3S_4)'-C_3S_4'C_5-S_3'C_4S_5'\right)\nonumber\\
     &&+S_5 S_6\, \left( (C S)_{1 4'}\, (C S)_{3 2'} +(C S)_{4 1'}\, (C S)_{2 3'}-2\right) \label{eq6.2}
    \end{eqnarray}
 \item[(d)] The point spectrum of the periodic quantum graph associated with trihexagonal tiling satisfies:\
     $\la=\rho^2$ lies in $\sig_p(H_{TH})$ if and only if $\rho$ satisfies
  \begin{eqnarray}
   0 &=&S_1 S_3 S_5+S_2 S_4 S_6+ (C S S S)_{1346} = S_1 S_4 S_5+S_2 S_3 S_6+(S_2 S_3 S_4 S_5)'\nonumber\\
    &=& S_1 S_2 S_4 +S_3 S_5 S_6+S_1 S_6\, (C S)_{2 5}+S_2 S_5\, (S_1S_6)' \label{eq6.3}
    \end{eqnarray}
   and
   \begin{eqnarray}
   0 &=&  -S_1 S_2 S_3-S_4 S_5 S_6 +(S_1 S_3 S_5+S_2 S_4 S_6)\, ((C S)_{16})'+(S_1 S_4 S_5+S_2 S_6 S_6)\, ((C S)_{25})'\nonumber\\
   && +(S_1 S_2 S_4+S_3 S_5 S_6)\, ((C S)_{34})'+ (S_1S_6)'\, (S_2S_5)'\, (C S)_{34}
   +(S_3S_4)'\, (C S)_{25}\, (C S)_{16}.\nonumber \\ \hfill\label{eq6.4}
   \end{eqnarray}
   \end{enumerate}
   \end{th6.1}
   \noindent
   {\bf Remark.} It is easy to verify that for the spectrum $ {\sigma(H_{trS})}$ related to truncated square tiling,
   $$
   \bigcup_{i,j,k,l,m\in \bfN_6} \{ \rho^2\in \bfR:\ S_i=S_j=S_k=S_l=S_m=0\}\ \subseteq \ \sig_p(H_{trS})\ .
   $$
   In particular, $\sig_p(H_{trS})$ contains those $\rho$ which satisfies $S_1=S_2=S_3=S_4=S_5=S_6=0$, whence the eigenfunction is known and there are infinitely many of them.
   Then for the spectrum $\sig_p(H_{TH})$ related to trihexagonal tiling,
   $$
   \bigcup_{i,j,k,l\in \bfN_6} \{ \rho^2\in \bfR:\ S_i=S_j=S_k=S_l=0\}\ \subseteq\   {\sigma_p(H_{TH})}\ .
   $$
   Certainly it includes the case when all $S_i$'s are zero.
   \begin{proof} (a) For the triangular tiling, let
   \begin{equation}
   \psi(\rho,\theta_1,\theta_2)=(S_1S_2S_3)'+(C S S)_{123}-2S_1S_2\, \cos(\theta_2-\theta_1)-2 S_2S_3\, \cos\theta_1-2 S_1S_3\, \cos\theta_2.
   \label{eq6.5}
   \end{equation}
   Now by the characterization of point spectrum, $\sig_p(H_{T})$ contains exactly those $\la=\rho^2$ which do not vary with $\Th$.
   Hence for any $\theta_1,\theta_2\in [-\pi,\pi]$,  {from Theorem~\ref{th2.1}},
   \begin{eqnarray*}
   0 &=& \frac{\partial\psi}{\partial\theta_1} = 2 S_1 S_2\, \sin(\theta_2-\theta_1)+2S_2S_3\, \sin\theta_1\\
   0 &=& \frac{\partial\psi}{\partial\theta_2} = -2 S_1 S_2\, \sin(\theta_2-\theta_1)+2S_2S_3\, \sin\theta_2
   \end{eqnarray*}
   Hence either $S_2=0$ or $S_1=S_3=0$.  In the first case, \eqref{eq6.5} implies
   $$
   S_1S_3\, (S_2'+C_2-2\cos\theta_2)=0,
   $$
   which means that either $S_1=0$ or $S_3=0$.  Thus part (a) is proved.
   \\[0.1in]
   (b)\ \ For the elongated triangular tiling, for any $\theta_1,\theta_2$,  {from Theorem~\ref{th3.1}},
   \begin{eqnarray}
   \frac{\partial\psi}{\partial\theta_1} &=& 2\left( (S_1S_2S_3S_4)'+(C S S S)_{1345}+S_1S_2S_5+S_1S_3S_4(C_2+S_5'+8\cos\theta_1)\right)\, \sin\theta_1\nonumber\\
    && +2 S_2S_3S_5\, \sin(\theta_1-\theta_2)\ = \ 0\label{eq6.6}\\
    \frac{\partial\psi}{\partial\theta_2} &=& 2 S_2S_4S_5\, \sin\theta_2-2S_2S_3S_5\, \sin(\theta_1-\theta_2)\ =\ 0.\label{eq6.7}
    \end{eqnarray}
    The equation \eqref{eq6.7} implies $S_2S_5=0$, or $S_3=S_4=0$, while \ $S_1S_3S_4=0$ \ follows from the first case.
    Eventually, we have that in any case there exists $i,j,k\in \bfN_5$ such that $S_i=S_j=S_k=0$.
    \\[0.1in]
    (c)\ \ For truncated square tiling, when $\theta_1,\theta_2\in [-\pi,\pi]$,  {from Theorem~\ref{th4.1} with $\displaystyle\psi=\Phi(\rho)\, \rme^{-i(\theta_1+\theta_2)}$},
    \begin{eqnarray*}
    \frac{\partial\psi}{\partial\theta_1} &=& 2\left( C_6S_1S_2+S_6(S_1S_2)'+(S_3S_4S_6)'\right)\, \sin\theta_1\\
    &&+2 S_2S_4\, \sin(\theta_1-\theta_2)+2S_1S_3\, \sin(\theta_1+\theta_2)\ = \ 0\\
    \frac{\partial\psi}{\partial\theta_2} &=& 2\left(S_5'S_2S_3+S_5(CS)_{23}+(CSS)_{145}\right)\, \sin\theta_2\\
    && -2S_2S_4\, \sin(\theta_1-\theta_2)+2S_1S_3\, \sin(\theta_1+\theta_2)\ = \ 0.
    \end{eqnarray*}
    It follows that \eqref{eq6.1} is valid.
    Together with Theorem \ref{th4.1}, \eqref{eq6.2} also follows.
    \\[0.1in]
    (d)\ \ For trihexagonal tiling,  {from Theorem~\ref{th5.1}},
    \begin{eqnarray*}
     \frac{\partial\psi}{\partial\theta_1} &=& -2\left( S_1S_2S_4+S_3S_5S_6+S_1S_6(CS)_{25}+S_2S_5(S_1S_6)'\right)\, \sin\theta_1\\
      &&-(S_1S_3S_5+S_2S_4S_6+(CSSS)_{1346})\, \sin(\theta_1-\theta_2)\ = \ 0\\
       \frac{\partial\psi}{\partial\theta_2} &=& -2(S_1S_4S_5+S_2S_3S_6+(S_2S_3S_4S_5)')\, \sin\theta_2\\
       && +(S_1S_3S_5+S_2S_4S_{ {6}}+(CSSS)_{1346})\, \sin(\theta_1-\theta_2)\ = \ 0
       \end{eqnarray*}
       These implies \eqref{eq6.3}, and together with Theorem \ref{th5.1} also yields \eqref{eq6.4}.
       \end{proof}

       The explicit dispersion relations given in Theorem \ref{th1.0}, assuming all $q_i$'s are identical and even, are simpler.
       They are all of the form $ {S^i\, p(S')=0}$, where $p$ is a polynomial  {whose coefficients depend on $\theta_1$ and $\theta_2$}.  Clearly, the point spectrum is defined by
       $S=0$, while the polynomials in $S'$ define the absolutely continuous spectra.
       \\[0.2in]
       \noindent
       {Proof of Theorem \ref{th1.1}:}\\
       For the triangular tiling,
       the  {dispersion relation \eqref{eq1.03}} implies $\rho^2\in \sig_{ac}(H_T)$ if and only if
       $$
       S'(a,\rho)=\frac{4}{3}\cos\frac{\theta_1}{2}\cos\frac{\theta_2}{2}\cos\frac{\theta_1-\theta_2}{2}-\frac{1}{3}
       =\frac{1}{6}\, \left| 1+\rme^{i\theta_1}+\rme^{i\theta_2}\right|^2-\frac{1}{2},
       $$
       by \eqref{eq1.2}. Since the range of $\displaystyle \left| 1+\rme^{i\theta_1}+\rme^{i\theta_2}\right|$ is  $[0,3]$. We conclude that
       $\rho^2\in \sig_{ac}(H_T)$ if and only if
       $$
       S'(a,\rho)\in [-\frac{1}{2},1].
       $$
       \hskip0.25in
       In the case of elongated triangular tiling,  {the dispersion relation \eqref{eq1.04}} is
       $$
       (S')^2-(\frac{4}{5}\cos\theta_1)\, S'-\frac{1}{25}\left(8\cos\frac{\theta_1}{2}\cos\frac{\theta_2}{2}\cos\frac{\theta_1-\theta_2}{2}
       -4\cos^2\theta_1+1\right)=0.
       $$
       It implies, through \eqref{eq1.2}, that $\rho^2\in\sig_{ac}(H_{eT})$ iff
       $$
       S'(a,\rho)=\frac{2}{5}\cos\theta_1\pm \frac{1}{5}\left| 1+\rme^{i\theta_1}+\rme^{i\theta_2}\right|\in [-\frac{3}{5},1].
       $$
       For truncated square tiling, if we let $\theta_1=\theta_2=\theta$, then  {\eqref{eq1.05} implies}
     \begin{eqnarray*}
     0&=& (9 S'^2-1)^2-(36 S'^2+24\cos\theta\, S'+4\, \cos^2\theta)\\
      &=& (9S'^2-1-6S'-2\cos\theta)(9S'^2-1+6S'+2\cos\theta).
      \end{eqnarray*}
     Thus,
     $$
     S'(a,\rho)=\frac{1}{3}({1\pm2\cos\frac{\theta}{2}})\in[\frac{1}{3},1]\cup[-\frac{1}{3},\frac{1}{3}]=[-\frac{1}{3},1];
     $$
     or
    $$
    S'(a,\rho)=\frac{1}{3}({-1\pm2\sin\frac{\theta}{2}})\in[-1,-\frac{1}{3}]\cup[-\frac{1}{3},\frac{1}{3}]=[-1,\frac{1}{3}].
    $$
    Also, if $|S'(a,\rho)|=1+\epsilon$ with $\epsilon>0$, then
   $$
   (9S'^2-1)^2-4(3S'+\cos\theta_1)(3S'+\cos\theta_2)\geq 192\ep+432 \ep^2+324\ep^3+81\ep^4>0.
   $$
   Hence we conclude that  $ \displaystyle  \sigma_{ac}(H_{trS})=\left\{S'(a,\rho)\in\left[-1,1\right]\right\}$.
   \par
  Finally the dispersion relation  {\eqref{eq1.06}} for trihexagonal tiling involves
     $$
      2S'^2-S'-\cos(\frac{\theta_1}{2})\cos(\frac{\theta_1}{2})\cos(\frac{\theta_1-\theta_2}{2})=0.
     $$
This implies that $\rho^2\in \sig_{ac}(H_{TH})$ iff
  \begin{eqnarray*}
 S'&=&\frac{1\pm\sqrt{1+8\cos\frac{\theta_1}{2}\cos\frac{\theta_2}{2}\cos\frac{\theta_1-\theta_2}{2}}}{4}\\
   &=&\frac{1}{4}(1\pm|1+{\rm e}^{i\theta_1}+{\rm e}^{i\theta_2}|)\\
  &\in&  { [\frac{1}{4},1]\cup[-\frac{1}{2},\frac{1}{4}]=[-\frac{1}{2},1]}.
  \end{eqnarray*}
  \begin{flushright}
  $\Box$
  \end{flushright}
       \hskip0.25in
       In case $q=0$,  {then $S'(a,\rho)=\cos(\rho a)$. Hence} from above, we have the following
\begin{enumerate}
\item[(a)] $\displaystyle \sigma_{ac}(H_T)=\sigma_{ac}(H_{TH})=\left[0,\frac{4\pi^2}{9a^2}\right]\cup\bigcup_{n=1}^\infty\left[(2n-\frac{2}{3})^2
    \frac{\pi^2}{a^2}, \, (2n+\frac{2}{3})^2\frac{\pi^2}{a^2}\right]$
\item[(b)] $\displaystyle \sigma_{ac}(H_{eT})=\left[0,\left(\frac{\eta}{a}\right)^2\right]
    \cup\bigcup_{n=1}^\infty\left[\left(\frac{2n\pi-\eta}{a}\right)^2
    ,\left(\frac{2n\pi+\eta}{a}\right)^2\right]$, $\eta=\cos^{-1}(-\frac{3}{5})\sim 2.2143$
\item[(c)]$\sigma_{ac}(H_{trS})=\mathbb{R}$
\end{enumerate}
\section{Concluding remarks}
  {Just like the case of trihexagonal tiling, we found infinitely many eigenfunctions for the point spectra for the periodic quantum graphs associated with
 the other three Archimedean tiling. Most eigenfunctions involves types (a) and (b) shown in Figure 3.  They are based on one triangle or one square of the tiling,
 and there are infinitely many of them.
  But for trihexagonal tiling we also have additional eigenfunctions based on a hexagon; and for truncated square tiling we have similar eigenfunctions based on a octagon.}

  {As a summary, we have proved that the dispersion relations for the periodic quantum graphs associated with the 4 Archimedean tilings are of the form
 $S^i (2S'-1)^j p(S')=0$ for some $i,j\geq 0$, and polynomial $p$.  Also their absolutely continuous spectra satisfy $\sig_{ac}\subset \phi^{-1}( [-1,1])$
 where $\phi(\rho)=S'(a,\rho)$. The spectra are all of a band and gap structure.
 In a next paper \cite{LJL}, we shall see that the situation for the other Archimedean tiling are essentially the same. Their dispersion relations are of
 the form  $S^i (2S'-1)^j(2S'+1)^k p(S')=0$, and their absolutely continuous spectra satisfy $\sig_{ac}\subset \phi^{-1}([-1,1])$.}

   {We would like to emphasize that the dispersion relations are derived with the help of Mathematica, because of the determinants of large matrices ($12\times
  12$ at the most). To guarantee that these relations are indeed correct, we also plug in arbitrary numerical values of $S_i(a,\rho)$, $S'_i(a,\rho)$ and
  $C_i(a,\rho)$ into the characteristic functions in Theorem~\ref{th2.1}-Theorem~\ref{th5.1}, to verify that they match with the determinants. Thus we are sure that our results are
  correct.}

  The dispersion relations obtained in this paper make it easy for the analysis of the spectrum.  We believe that there will be a lot of interesting properties
 about these quantum graphs to be dug out.  We also emphasize that our method, through some hard work,  may allow $q_i$'s to be neither identical, nor even, as given in
 Theorems \ref{th2.1}=\ref{th5.1}. For the sake of completeness, we also include the cases of square tiling \cite{LLW},
 \begin{eqnarray*}
 \lefteqn{ (C_1S_2S_3S_4)'+(S_1C_2S_3S_4)'+(S_1S_2C_3S_4)'+(S_1S_2S_3C_4)'}\\
 &=&2[(S_1S_3+S_2S_4)\cos\theta_1+(S_1S_2+S_3S_4)\cos\theta_2+S_2S_3\cos(\theta_1+\theta_2)+S_1S_4\cos(\theta_1-\theta_2)].
 \end{eqnarray*}
 \vskip0.2in
   \section*{Acknowledgements}
 {The authors are grateful to the anonymous reviewers for their very careful reading and constructive comments. Their suggestions
 help us to improve the paper significantly and correct some mistakes.
We also thank Vyacheslav Pivovarchik and Min-Jei Huang for stimulating discussions.
The authors are partially supported by Ministry of Science and Technology, Taiwan, under contract number MOST105-2115-M-110-004.}
  \newpage
 \begin{table}[h!]
 \centering
 \caption{Table of 11 Archimedean tilings}
 $\begin{array}{|c||c|c|c|}
 \hline
 \text{Name} & \text{Triangular tiling} & \text{{\small Snub trihexagonal tiling}} & \text{{\small Elongated triangular tiling}} \\ \hline
 \text{Notation} & \left(3^6\right) & \left(3^4,6\right) & \left(3^3,4^2\right) \\ \hline
 \  & \includegraphics[width=3cm]{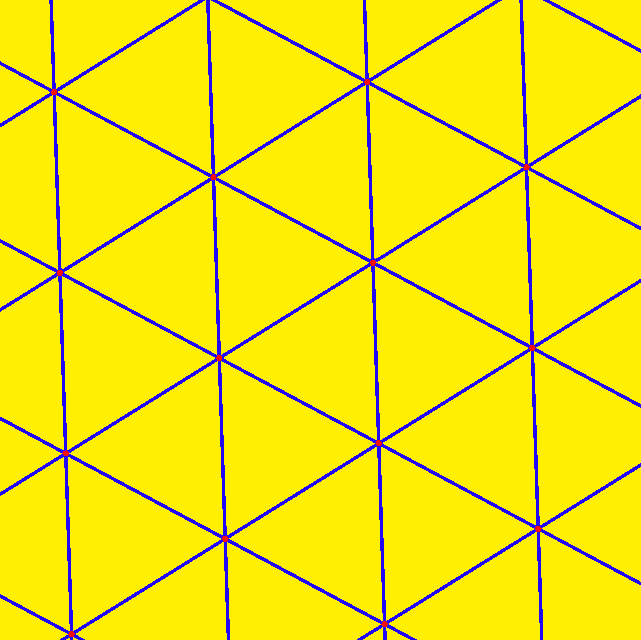} & \includegraphics[width=3cm]{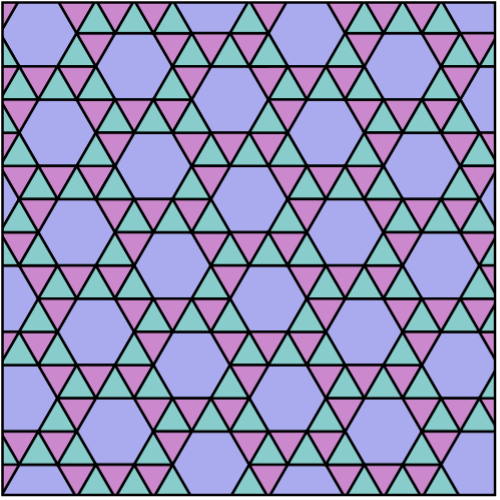} & \includegraphics[width=3cm]{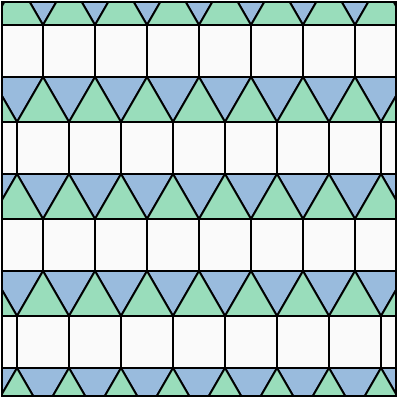}  \\ \hline \hline
 \text{Name} & \text{Snub square tiling} & \text{Trihexagonal tiling} & \text{{\small Rhombi-trihexagonal tiling}} \\ \hline
 \text{Notation} & \left(3^2,4,3,4\right) & (3, 6, 3, 6) & (3,4,6,4) \\ \hline
 \  & \includegraphics[width=3cm]{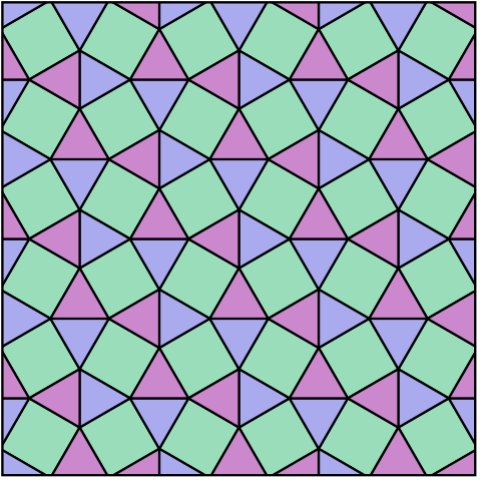} & \includegraphics[width=3cm]{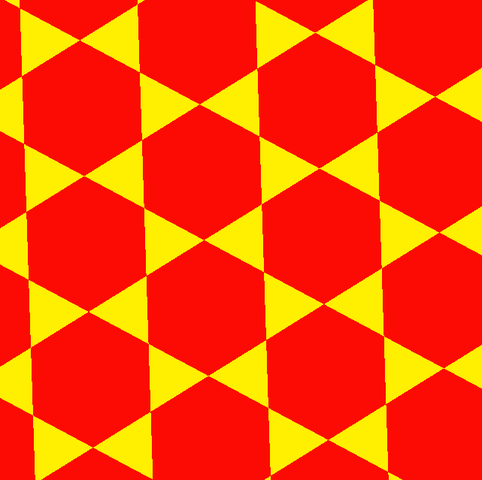} & \includegraphics[width=3cm]{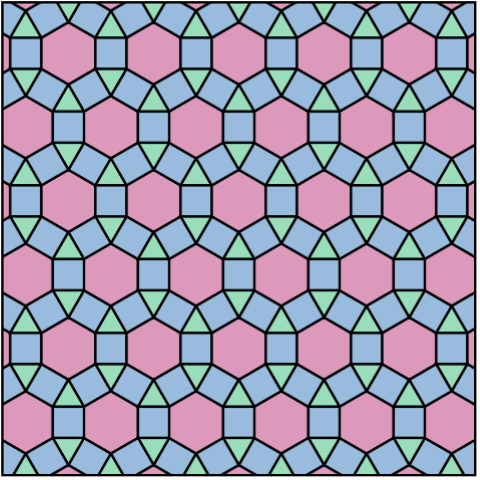} \\ \hline \hline
 \text{Name} & \text{{\small Truncated hexagonal tiling}} & \text{Square tiling} & \text{{\small Truncated trihexagonal tiling}} \\ \hline
 \text{Notation} & \left(3,12^2\right) & \left(4^4\right) & (4,6,12) \\ \hline
 \  & \includegraphics[width=3cm]{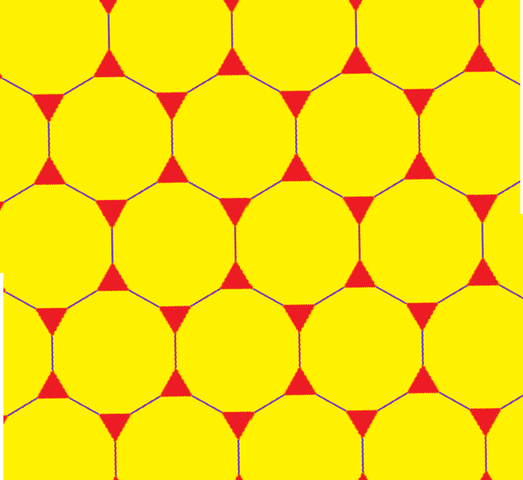} & \includegraphics[width=3cm]{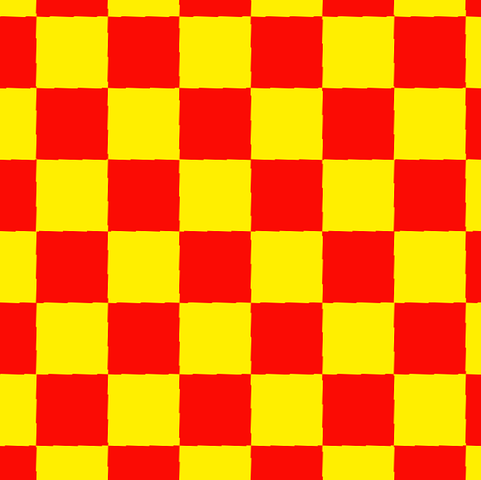} & \includegraphics[width=3cm]{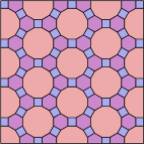} \\ \hline \hline
 \text{Name} & \text{{\small Truncated square tiling}} & \text{Hexagonal tiling} & \ \\ \hline
 \text{Notation} & \left(4,8^2\right) & (6^3) & \ \\ \hline
 \  & \includegraphics[width=3cm]{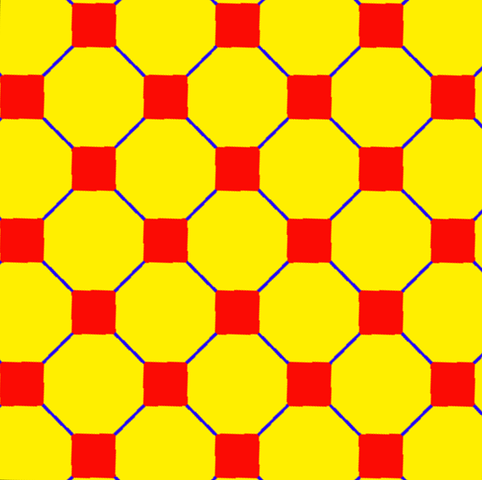} & \includegraphics[width=3cm]{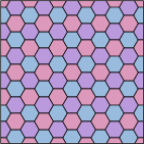} & \ \\ \hline
 \end{array}\label{Tab5.1}$
 \end{table}

 \end{document}